\documentclass[a4paper,10pt,twoside]{article}

\usepackage{amssymb}
\usepackage{amsthm}
\usepackage{amsmath}
\usepackage{enumerate}
\usepackage{hyperref}

\newcommand{\Pb}{\mathbf{P} \left\{ }
\newcommand{\lb}{\left\{ }
\newcommand{\rb}{\right\} }
\newcommand{\lpa}{\left( }
\newcommand{\rpa}{\right) }
\newcommand{\lsb}{\left[ }
\newcommand{\rsb}{\right] }
\newcommand{\la}{\left| }
\newcommand{\ra}{\right| }
\newcommand{\E}{\mathbf{E}}
\newcommand{\Var}{\mathbf{Var}}

\newcommand{\Hh}{H-\half}
\newcommand{\Hth}{H-\threehalf}
\newcommand{\Ga}{\frac{1}{\Gamma \left(H+\half\right)}}
\newcommand{\Gaa}{\Gamma\left(H+\half\right)}
\newcommand{\WH}{W^{(H)}}
\newcommand{\WHd}{W_{\delta}^{(H)}}
\newcommand{\BmH}{B_m^{(H)}}
\newcommand{\BmpH}{B_{m+1}^{(H)}}
\newcommand{\BmHd}{B_{m,\delta}^{(H)}}
\newcommand{\trj}{t_{r+\frac j4}}
\newcommand{\tm}{t_{(m)}}
\newcommand{\dm}{\delta_{(m)}}
\newcommand{\dmi}{\left(\frac{-1}{\delta}\right)_{(m)}}
\newcommand{\dmit}{\textstyle \left(\frac{-1}{\delta}\right)_{(m)}}
\newcommand{\quart}{\frac14 }
\newcommand{\half}{\frac12 }
\newcommand{\threehalf}{\frac32 }

\newcommand*{\di}{\, \mathrm{d} }

\newtheorem{thmn}{Theorem}

\newtheorem{lem}{Lemma}

\begin{document}


\title{Strong Approximation of Fractional Brownian Motion by Moving
Averages of Simple Random Walks}
\author{Tam\'as Szabados \\ Department of Mathematics,
 Technical University of Budapest \\
 Egry u 20-22, H \'ep. V em.,
 Budapest, 1521 \\
 Hungary\\
 E-mail: szabados@math.bme.hu}
\date{}


\maketitle

\begin{center} \emph{Dedicated to P\'al R\'ev\'esz on the occasion of
his 65th birthday }
\end{center}


\begin{abstract}
The fractional Brownian motion is a generalization of ordinary Brownian
motion, used particularly when long-range dependence is required. Its
explicit introduction is due to B.B. Mandelbrot and J.W. van Ness (1968)
as a self-similar Gaussian process $\WH (t)$ with stationary increments.
Here self-similarity means that $(a^{-H}\WH(at): t \ge 0)
\stackrel{d}{=} (\WH(t): t \ge 0)$, where $H\in (0, 1)$ is the Hurst
parameter of fractional Brownian motion.

F.B. Knight gave a construction of ordinary Brownian motion as a limit
of simple random walks in 1961. Later his method was simplified by
P. R\'ev\'esz (1990) and then by the present author (1996). This
approach is quite natural and elementary, and as such, can be extended
to more general situations. Based on this, here we use moving averages
of a suitable nested sequence of simple random walks that almost surely
uniformly converge to fractional Brownian motion on compacts when $H \in
(\quart , 1)$. The rate of convergence proved in this case is
$O(N^{-\min(H-\quart,\quart)}\log N)$, where $N$ is the number of steps
used for the approximation.

If the more accurate (but also more intricate) Koml\'os, Major,
Tusn\'ady (1975, 1976) approximation is used instead to embed random
walks into ordinary Brownian motion, then the same type of moving
averages almost surely uniformly converge to fractional Brownian motion
on compacts for any $H \in (0, 1)$. Moreover, the convergence rate is
conjectured to be the best possible $O(N^{-H}\log N)$, though only
$O(N^{-\min(H,\half)}\log N)$ is proved here.

\end{abstract}






\section{Fractional Brownian motion}

The fractional Brownian motion (fBM) is a generalization of ordinary
Brownian motion (BM) used particularly when long-range dependence is
essential. Though the history of fBM can be traced back to \cite{Kol40}
and others, its explicit introduction is due to \cite{Man68}.
Their intention was to define a \emph{self-similar},
centered Gaussian process $\WH (t)$ $(t\ge 0)$ with stationary but not
independent increments and with continuous sample paths a.s. Here
self-similarity means that for any $a>0$,
\begin{equation}
\lpa a^{-H} \WH (at) : t\ge 0 \rpa \stackrel{d}{=} \lpa \WH (t) :
t\ge 0 \rpa , \label{eq:ss}
\end{equation}
where $H \in (0,1)$ is the \emph{Hurst parameter} of the fBM and
$\stackrel{d}{=}$ denotes equality in distribution. They showed that
these properties characterize fBM. The case $H=\half$ reduces to
ordinary BM with independent increments, while the cases $H<\half$ and
$H>\half$ give negatively, respectively, positively correlated
increments, see \cite{Man68}. It seems that in the
applications of fBM the case $H>\half$ is the most frequently used.

\cite{Man68} gave the following explicit
representation of fBM as a moving average of ordinary, but two-sided BM
$W(s),~s\in \mathbb{R}$:
\begin{equation}
\WH (t) = \Ga \int_{-\infty}^t \lsb (t-s)^{\Hh}
- (-s)_+^{\Hh} \rsb \, \di W(s) \qquad (t \ge 0), \label{eq:fbrow}
\end{equation}
where $(x)_+ = \max(x,0)$.

The idea of (\ref{eq:fbrow}) is related to \emph{deterministic
fractional calculus}, which has an even longer history than fBM, going
back to Liouville, Riemann, and others, see in \cite{Sam93}. Its
simplest case is when a continuous function $f$ and a positive integer
$\alpha $ are given. Then an induction with integration by parts can
show that
\[ f_{\alpha}(t) = \frac{1}{\Gamma(\alpha)} \int_0^t (t-s)^{\alpha -1}
 f(s) \di s \]
is the order $\alpha $ iterated antiderivative (or order $\alpha $
integral) of $f$. On the other hand, this integral is well-defined for
non-integer positive values of $\alpha$ as well, in which case it can be
called a fractional integral of $f$.

So, heuristically, the main part of (\ref{eq:fbrow}),
\[ W_{\alpha}(t) = \frac{1}{\Gamma(\alpha)} \int_0^t (t-s)^{\alpha - 1}
W'(s) \di s = \frac{1}{\Gamma(\alpha)} \int_0^t (t-s)^{\alpha -1} \di W(s) \]
is the order $\alpha$ integral of the (in ordinary sense non-existing)
white noise process $W'(t)$. Thus the fBM $\WH (t)$ can be considered as
a stationary-increment modification of the fractional integral
$W_{\alpha}(t)$ of the white noise process, where $\alpha = H+\half \in
(\half,\threehalf)$.


\section{Random walk construction of ordinary Brownian motion}

It is interesting that a very natural and elementary construction of
ordinary BM as a limit of random walks (RWs) appeared relatively
late. The mathematical theory of BM began around 1900 with the works of
Bachelier, Einstein, Smoluchowski, and others. The first existence
construction was given by Wiener (1921, 1923) that was followed by
several others later. F.B. Knight (1961) introduced the first
construction by random walks that was later simplified by P.  R\'ev\'esz
(1990). The present author was fortunate enough to hear this version of
the construction directly from P\'al R\'ev\'esz in a seminar at the
Technical University of Budapest a couple of years before the
publication of R\'ev\'esz's book in 1990 and got immediately fascinated
by it. The result of an effort to further simplify it appeared in
\cite{Sza96}.
From now on, the expression \emph{RW construction} will always refer to
the version discussed in the latter. It is asymptotically
equivalent to applying \cite{Sko65} embedding to find a nested
dyadic sequence of RWs in BM, see Theorem~4 in \cite{Sza96}. As
such, it has some advantages and disadvantages compared to
the celebrated best possible approximation by BM of partial sums of
random variables with moment generator function finite around the
origin. The latter was obtained by Koml\'os, Major, and Tusn\'ady (1975,
1976), and will be abbreviated \emph{KMT approximation} in the sequel.
The main advantages of the RW construction are that it is elementary,
explicit, uses only past values to construct new ones, easy to implement
in practice, and very suitable for approximating stochastic integrals,
see Theorem~6 in \cite{Sza96} and also \cite{Sza90}.
Recall that the KMT approximation constructs partial sums (e.g. a simple
symmetric RW) from BM itself (or from an i.i.d. sequence of standard
normal random variables) by an intricate sequence of conditional
quantile transformations. To construct any new value it uses to whole
sequence (past and future values as well). On the other hand, the major
weakness of the RW construction is that it gives a rate of convergence
$O(N^{-\quart}\log N)$, while the rate of the KMT approximation is the
best possible $O(N^{-\half} \log N)$, where $N$ is the number of steps
(terms) considered in the RW.

In the sequel first the main properties of the above-mentioned RW
construction are summarized. Then this RW construction is used to define
an approximation, similar to (\ref{eq:fbrow}), of fBM by moving averages
of the RW. The convergence and the error of this approximation are
discussed next.  As a consequence of the relatively weaker approximation
properties of the RW construction, the convergence to fBM will be
established only for $H \in (\quart, 1)$, and the rate of
convergence will not be the best possible either. To compensate for
this, at the end of the paper we discuss the convergence and error
properties of a similar construction of fBM that uses the KMT
approximation instead, which converges for all $H \in (0, 1)$ and whose
convergence rate can be conjectured to be the best possible when
approximating fBM by moving averages of RWs.

The RW construction of BM summarized here is taken from
\cite{Sza96}.
We start with an infinite matrix of i.i.d. random variables $X_m(k)$,
\[ \Pb X_m(k) = 1 \rb = \Pb X_m(k) = -1 \rb =
\half \qquad (m\ge 0, k\ge 1), \]
defined on the same underlying probability space
$(\Omega, \mathcal{A}, \mathbf{P})$. Each row of this matrix is a basis
of an approximation of BM with a certain dyadic step size $\Delta t =
2^{-2m}$ in time and a corresponding step size $\Delta x = 2^{-m}$ in
space, illustrated by the next table.

\begin{table}[h]
\caption{The starting setting for the RW construction of BM}
\[
\begin{array}{|c|c|l|l|}
\hline
\Delta t &\Delta x & \mbox{i.i.d. sequence} & \mbox{RW} \\
\hline
1&1& X_0(1), X_0(2), X_0(3), \dots & S_0(n) = \sum_{k=1}^n X_0(k) \\
2^{-2}&2^{-1}& X_1(1), X_1(2), X_1(3), \dots & S_1(n) = \sum_{k=1}^n
X_1(k) \\
2^{-4}&2^{-2}& X_2(1), X_2(2), X_2(3), \dots & S_2(n) = \sum_{k=1}^n
X_2(k) \\
\vdots & \vdots & \vdots & \vdots \\
\hline
\end{array}
\]
\end{table}

The second step of the construction is \emph{twisting}. From the
independent random walks (i.e. from the rows of Table 1), we want to
create dependent ones so that after shrinking temporal and spatial step
sizes, each consecutive RW becomes a refinement of the previous one.
Since the spatial unit will be halved at each consecutive row, we define
stopping times by $T_m(0)=0$, and for $k \ge 0$,
\[ T_m(k+1) = \min \lb n : n > T_m(k), |S_m(n) - S_m(T_m(k))| = 2 \rb
\qquad (m \ge 1). \]
These are the random time instants when a RW visits even integers,
different from the previous one. After shrinking the spatial unit by
half, a suitable modification of this RW will visit the same integers in
the same order as the previous RW. (This is what we call a refinement.)
We will operate here on each point $\omega \in \Omega $ of the sample
space separately, i.e. we fix a sample path of each RW appearing in
Table 1. Thus each \emph{bridge} $S_m(T_m(k+1)) - S_m(T_m(k))$ has to
mimic the corresponding step $X_{m-1}(k+1)$ of the previous RW. We
define twisted RWs $\tilde S_m$ recursively for $m=1,2,3,\dots $ using
$\tilde S_{m-1}$, starting with $\tilde S_0(n) = S_0(n)$ $(n \ge 0)$.
With each fixed $m$ we proceed for $k=0,1,2,\dots $ successively, and
for every $n$ in the corresponding bridge, $T_m(k) < n \le T_m(k+1)$.
Any bridge is flipped if its sign differs from the desired (Figs. 1-3):
\[ \tilde X_m(n) = \left\{
\begin{array}{rl}
X_m(n)& \mbox{ if } \ S_m(T_m(k+1)) - S_m(T_m(k))
 = 2 \tilde X_{m-1}(k+1), \\
-X_m(n)& \mbox{ otherwise,}
\end{array} \right. \]
and then $\tilde S_m(n) = \tilde S_m(n-1) + \tilde X_m(n)$.  Then each
$\tilde S_m(n)$ $( n \ge 0)$ is still a simple, symmetric RW, see
Lemma~1 in \cite{Sza96}. Moreover, the twisted RWs have the
desired refinement property:
\begin{equation}
\half \tilde S_m(T_m(k)) = \tilde S_{m-1}(k) \qquad
(m \ge 1, k \ge 0). \label{eq:refis}
\end{equation}


\begin{figure}[!ht]
\begin{picture}(50,40)(-10,-10)
\thinlines
\put (-5,0){\vector(1,0){40}}
\put (0,-5){\vector(0,1){30}}
\put (10,0){\line(0,1){2}}
\put (20,0){\line(0,1){2}}
\put (30,0){\line(0,1){2}}
\put (0,10){\line(1,0){2}}
\put (0,20){\line(1,0){2}}
\thicklines
\put (40,0){\makebox(0,0){$t$}}
\put (0,30){\makebox(0,0){$x$}}
\put (10,-3){\makebox(0,0){1}}
\put (20,-3){\makebox(0,0){2}}
\put (30,-3){\makebox(0,0){3}}
\put (-2,10){\makebox(0,0){1}}
\put (-2,20){\makebox(0,0){2}}
\put (0,0){\line(1,1){20}}
\put (20,20){\line(1,-1){10}}
\put (0,0){\circle*{3}}
\put (10,10){\circle*{3}}
\put (20,20){\circle*{3}}
\put (30,10){\circle*{3}}
\end{picture}
\caption{$B_0(t;\omega ) = S_0(t;\omega )$.}
\end{figure}
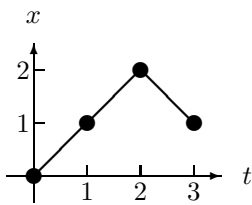

\begin{figure}[!ht]
\begin{picture}(175,55)(-5,-25)
\thinlines
\put (-5,0){\vector(1,0){170}}
\put (0,-25){\vector(0,1){50}}
\put (10,0){\line(0,1){2}}
\put (0,20){\line(1,0){2}}
\put (0,-20){\line(1,0){2}}
\thicklines
\put (170,0){\makebox(0,0){$t$}}
\put (0,30){\makebox(0,0){$x$}}
\put (10,-3){\makebox(0,0){1}}
\put (-2,20){\makebox(0,0){2}}
\put (-3,-20){\makebox(0,0){-2}}
\put (0,0){\line(1,-1){10}}
\put (10,-10){\line(1,1){10}}
\put (20,0){\line(1,-1){10}}
\put (30,-10){\line(1,1){20}}
\put (50,10){\line(1,-1){30}}
\put (80,-20){\line(1,1){20}}
\put (100,0){\line(1,-1){10}}
\put (110,-10){\line(1,1){20}}
\put (130,10){\line(1,-1){10}}
\put (140,0){\line(1,1){20}}
\put (0,0){\circle*{3}}
\put (80,-20){\circle*{3}}
\put (100,0){\circle*{3}}
\put (160,20){\circle*{3}}
\end{picture}
\caption{$S_1(t;\omega )$.}
\end{figure}
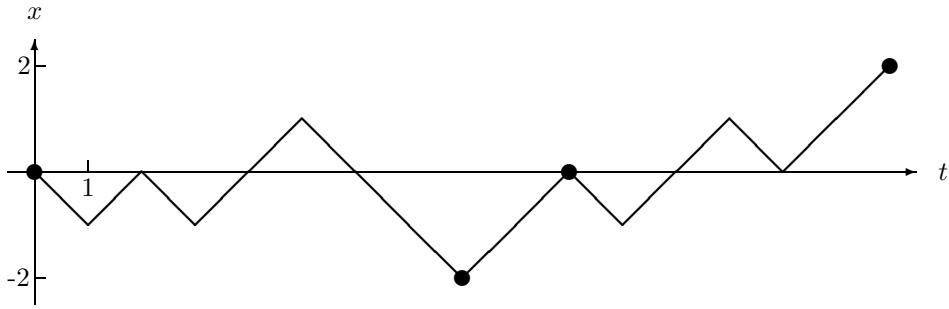

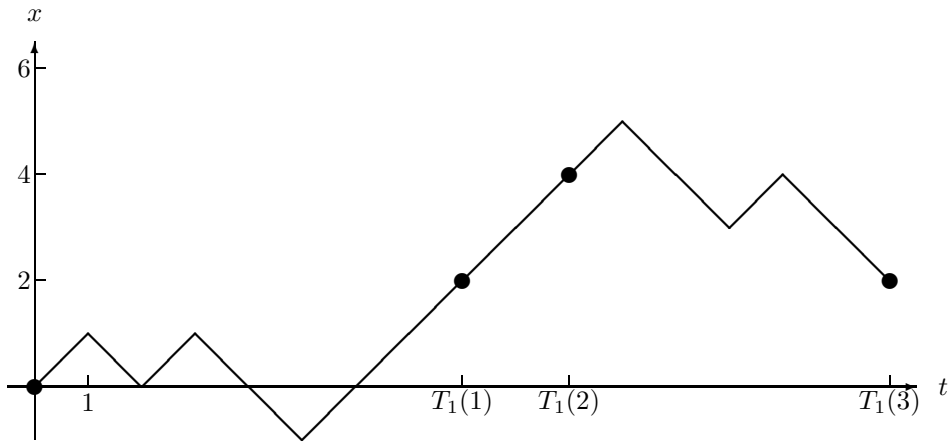
\begin{figure}[!ht]
\begin{picture}(175,80)(-5,-10)
\thinlines
\put (-5,0){\vector(1,0){170}}
\put (0,-10){\vector(0,1){75}}
\put (10,0){\line(0,1){2}}
\put (80,0){\line(0,1){2}}
\put (100,0){\line(0,1){2}}
\put (160,0){\line(0,1){2}}
\put (0,20){\line(1,0){2}}
\put (0,40){\line(1,0){2}}
\put (0,60){\line(1,0){2}}
\thicklines
\put (170,0){\makebox(0,0){$t$}}
\put (0,70){\makebox(0,0){$x$}}
\put (10,-3){\makebox(0,0){1}}
\put (80,-3){\makebox(0,0){$T_1(1)$}}
\put (100,-3){\makebox(0,0){$T_1(2)$}}
\put (160,-3){\makebox(0,0){$T_1(3)$}}
\put (-2,20){\makebox(0,0){2}}
\put (-2,40){\makebox(0,0){4}}
\put (-2,60){\makebox(0,0){6}}
\put (0,0){\line(1,1){10}}
\put (10,10){\line(1,-1){10}}
\put (20,0){\line(1,1){10}}
\put (30,10){\line(1,-1){20}}
\put (50,-10){\line(1,1){60}}
\put (110,50){\line(1,-1){20}}
\put (130,30){\line(1,1){10}}
\put (140,40){\line(1,-1){20}}
\put (0,0){\circle*{3}}
\put (80,20){\circle*{3}}
\put (100,40){\circle*{3}}
\put (160,20){\circle*{3}}
\end{picture}
\caption{$\tilde S_1(t;\omega )$.}
\end{figure}

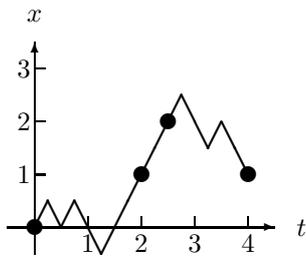
\begin{figure}[!ht]
\begin{picture}(60,50)(-10,-10)
\thinlines
\put (-5,0){\vector(1,0){50}}
\put (0,-5){\vector(0,1){40}}
\put (10,0){\line(0,1){2}}
\put (20,0){\line(0,1){2}}
\put (30,0){\line(0,1){2}}
\put (40,0){\line(0,1){2}}
\put (0,10){\line(1,0){2}}
\put (0,20){\line(1,0){2}}
\put (0,30){\line(1,0){2}}
\thicklines
\put (50,0){\makebox(0,0){$t$}}
\put (0,40){\makebox(0,0){$x$}}
\put (10,-3){\makebox(0,0){1}}
\put (20,-3){\makebox(0,0){2}}
\put (30,-3){\makebox(0,0){3}}
\put (40,-3){\makebox(0,0){4}}
\put (-2,10){\makebox(0,0){1}}
\put (-2,20){\makebox(0,0){2}}
\put (-2,30){\makebox(0,0){3}}
\put (0,0){\line(1,2){2.5}}
\put (2.5,5){\line(1,-2){2.5}}
\put (5,0){\line(1,2){2.5}}
\put (7.5,5){\line(1,-2){5}}
\put (12.5,-5){\line(1,2){15}}
\put (27.5,25){\line(1,-2){5}}
\put (32.5,15){\line(1,2){2.5}}
\put (35,20){\line(1,-2){5}}
\put (0,0){\circle*{3}}
\put (20,10){\circle*{3}}
\put (25,20){\circle*{3}}
\put (40,10){\circle*{3}}
\end{picture}
\caption{$B_1(t;\omega )$.}
\end{figure}


The last step of the RW construction is \emph{shrinking}. The sample
paths of $\tilde S_m(n)$ $(n \ge 0)$ can be extended to continuous
functions by linear interpolation, this way one gets $\tilde S_m(t)$
$(t \ge 0)$ for real $t$. Then we define the \emph{$m$th approximation
of BM} (see Fig. 4) by
\begin{equation}
B_m(t) = 2^{-m} \tilde S_m(t 2^{2m}). \label{eq:shrink}
\end{equation}

Compare three steps of a sample path of the first approximation
$B_0(t;\omega )$ and the corresponding part of the second approximation
$B_1(t; \omega )$ on Figs. 1 and 4. The second visits the same integers
(different from the previous one) in the same order as the first, so
mimics the first, but the corresponding time instants differ in general:
$ 2^{-2} T_1(k) \ne k$. Similarly, (\ref{eq:refis}) implies the general
\emph{refinement property}
\begin{equation}
B_{m+1} \lpa T_{m+1}(k) \, 2^{-2(m+1)} \rpa = B_m \lpa k 2^{-2m}
\rpa \qquad (m \ge 0, k \ge 0), \label{eq:refib}
\end{equation}
but there is a \emph{time lag}
\begin{equation}
T_{m+1}(k) \  2^{-2(m+1)} - k 2^{-2m} \ne 0 \label{eq:tlag}
\end{equation}
in general.

The basic idea of the RW construction of BM is that these time lags
become uniformly small if $m$ gets large enough. It can be proved by the
following simple lemma.

\begin{lem} \label{le:basic}
Suppose that $X_1, X_2, \dots , X_N$ is an i.i.d. sequence of random
variables, $\E (X_k)=0$, $\Var (X_k)=1$, and their moment generating
function $\E \lpa e^{u X_k} \rpa$ is finite for $|u| \le u_0, \ u_0
> 0$. Let $S_j = X_1+ \cdots +X_j, \ 1 \le j \le N$. Then for any $C>1$
and $N \ge N_0(C)$ one has
\[ \Pb \max_{0 \le j \le N} |S_j| \ge (2C N \log N)^{\half} \rb \le
2 N^{1-C}. \]
\end{lem}

This basic fact follows from a large deviation inequality, see e.g.
Section XVI,6 in \cite{Fel66}. Lemma~\ref{le:basic} easily implies
the uniform smallness of time lags in (\ref{eq:tlag}).

\begin{lem} \label{le:timelag}
For any $K > 0$, $C>1$, and for any $m \ge m_0(C)$, we have
\begin{multline*}
\Pb \max_{0 \le k 2^{-2m} \le K} |T_{m+1}(k) \, 2^{-2(m+1)} -
k 2^{-2m}| \ge \lpa \threehalf C K \log_* K \rpa ^{\half}
m^{\half} \, 2^{-m} \rb \\
\le 2 \, (K 2^{2m})^{1-C},
\end{multline*}
where $\log_*(x) = \max (1, \log x)$.
\end{lem}

Not surprisingly, this and the refinement property (\ref{eq:refib}) imply
the uniform closeness of two consecutive approximations of BM if $m$
is large enough.

\begin{lem} \label{le:apprb}
For any $K > 0$, $C > 1$, and $m \ge m_1(C)$, we have
\begin{multline*}
\Pb \max_{0 \le k 2^{-2m} \le K} |B_{m+1}(k 2^{-2m}) -
 B_m(k 2^{-2m})| \ge K^{\quart} (\log_*K)^{\frac{3}{4}} m 2^{-\frac{m}{2}}
\rb \\
\le 3 (K 2^{2m})^{1-C}.
\end{multline*}
\end{lem}

This lemma ensures the a.s. uniform convergence of the RW
approximations on compact intervals and it is clear that the limit
process is the Wiener process (BM) with continuous sample
paths almost surely.

\begin{thmn} \label{th:wiener}
The RW approximation $B_m(t)$  $(t \ge 0, m=0,1,2, \dots)$ a.s.
uniformly converges to a Wiener process $W(t)$  $(t \ge 0)$ on any
compact interval $[0, K], K > 0$.  For any $K > 0$, $C \ge 3/2$, and
for any $m \ge m_2(C)$, we have
\[ \Pb \max_{0 \le t \le K} |W(t) - B_m(t)| \ge
 K^{\quart} (\log_*K)^{\frac{3}{4}} m 2^{-\frac{m}{2}} \rb
\le 6 (K 2^{2m})^{1-C}. \]
\end{thmn}

The results quoted above correspond to Lemmas 2-4 and Theorem 3 in
\cite{Sza96}. We mention that the statements presented here are given
in somewhat sharper forms, but they can be read easily from the proofs
in the above reference.


\section{A pathwise approximation of fractional Brownian motion}

An almost surely convergent pathwise construction of fBM was given by
Carmona and Coutin (1998) representing fBM as a linear functional of an
infinite dimensional Gaussian process. Another pathwise construction
was given by Decreusefond and \"Ust\"unel (1998, 1999) which
converges in the $L^2$ sense. This construction uses discrete
approximations of the moving average representation of fBM
(\ref{eq:fbrow}), based on deterministic partitions of the time axis.
More exactly, (\ref{eq:fbrow}) is substituted by an integral over the
compact interval $[0, t]$, but with a more complicated kernel containing
a hypergeometric function too.

The approximation of fBM discussed here will also be a discrete version
of the moving average representation (\ref{eq:fbrow}) of fBM, but dyadic
partitions are taken on the spatial axis of BM and so one gets random
partitions on the time axis. This is asymptotically a Skorohod-type
embedding of nested RWs into BM. As a result, instead of integral we
have sum, and BM is substituted by the nested, refining sequence of its
RW approximations discussed in the previous section.  Since
(\ref{eq:fbrow}) contains two-sided BM, we need two such sequences: one
for the right and one for the left half-axis. From now on, we are going
to use the following notations: $m \ge 0$ is an integer, $\Delta t =
2^{-2m}$, $t_x = x \Delta t$ $(x \in \mathbb{R})$.
Then by definition, the \emph{$m$th approximation of fBM} is:
$\BmH(0)=0$, and for positive integers $k$,
\begin{multline}
\BmH (t_k) \\
= \Ga \sum_{r=-\infty }^{k-1} \lsb (t_k - t_r)^{\Hh}
-(-t_r)_+^{\Hh} \rsb \lsb B_m(t_r+\Delta t) - B_m(t_r) \rsb \\
= \frac{2^{-2Hm}}{\Gaa} \sum_{r=-\infty }^{k-1} \lsb (k - r)^{\Hh}
-(-r)_+^{\Hh} \rsb \tilde X_m(r+1),  \label{eq:afbrow1}
\end{multline}
where the convention $0^{\Hh}=0$ is applied even for negative exponents.

$\BmH $ is \emph{well-defined}, since the ``infinite part''
\[ \sum_{r=-\infty }^{-1} \lsb (k - r)^{\Hh} - (-r)^{\Hh} \rsb
\tilde X_m(r+1) =: \sum_{r=-\infty }^{-1} Y_{k,-r} \]
converges a.s. to a random variable $Z_k$ by Kolmogorov's ``three-series
theorem'': $\E(Y_{k,v}) = 0$ and
\begin{multline*}
\sum_{v=1 }^{\infty} \Var (Y_{k,v}) = \sum_{v=1}
^{\infty } v^{2H-1} \lsb \lpa 1 + \frac {k}{v} \rpa ^{\Hh} - 1 \rsb^2
\sim \sum_{v=1}^{\infty } \frac {\mbox {const}}{v^{3 - 2H}} < \infty . 
\end{multline*}

It is useful to write $\BmH $ in another form applying a discrete
version of integration by parts. Starting with (\ref{eq:afbrow1})
and rearranging it according to $B_m(t_r)$, one obtains for $k \ge 1$
that
\begin{equation}
\BmH (t_k) = \sum_{r=-\infty }^k  \frac{h(t_r - \Delta t, t_k) -
h(t_r, t_k)}{\Delta t} \, B_m(t_r) \, \Delta t , \label{eq:afbrow2}
\end{equation}
where we introduced the kernel
\begin{equation}
h(s,t) = \Ga \lsb (t-s)^{\Hh} - (-s)_+^{\Hh} \rsb \qquad (s \le t).
\label{eq:hst}
\end{equation}
This way we have got a discrete version of
\begin{equation}
\WH (t) = \frac{-1}{\Gaa} \int_{-\infty}^t \frac{\di}{\di s}
\lsb (t-s)^{\Hh} - (-s)_+^{\Hh} \rsb W(s) \, \di s, \label{eq:fbrow2}
\end{equation}
which is what one obtains from (\ref{eq:fbrow}) using a formal
integration by parts (cf. Lemma~\ref{le:Gauss} below).

To support the above definition we show that $\BmH$ has properties
analogous to the characterizing properties of fBM in a discrete setting.

\emph{(a)} $\BmH $ is \emph{centered} (clear from its definition) and has
\emph{stationary increments}: if $k_0$ and $k$ are non-negative integers,
then (substituting $u=r-k_0$)
\begin{eqnarray*}
\lefteqn{\BmH (t_{k_0} + t_k) - \BmH (t_{k_0})} \\
&=& \frac{2^{-2Hm}}{\Gaa} \lb \sum_{r=0}^{k_0+k-1} (k_0+k-r)^{\Hh}
\tilde X_m(r+1) \right. \\
&&- \left. \sum_{r=0}^{k_0} (k_0-r)^{\Hh} \tilde X_m(r+1)
\right. \\
&&+ \left. \sum_{r=-\infty}^{-1} \lsb (k_0+k-r)^{\Hh} - (k_0 - r)^{\Hh}
\rsb \tilde X_m(r+1) \rb \\
&=& \frac{2^{-2Hm}}{\Gaa} \lb \sum_{u=-k_0}^{k-1} (k-u)^{\Hh}
\tilde X_m(k_0+u+1) \right. \\
&&- \left. \sum_{u=-k_0}^{0} (-u)^{\Hh} \tilde X_m(k_0+u+1)
\right. \\
&& + \left. \sum_{u=-\infty}^{-k_0-1} \lsb (k-u)^{\Hh} - (-u)^{\Hh} \rsb
\tilde X_m(k_0+u+1) \rb \\
&=& \frac{2^{-2Hm}}{\Gaa} \sum_{u=-\infty}^{k-1} \lsb (k-u)^{\Hh}
- (-u)_+^{\Hh}\rsb \tilde X_m(k_0+u+1) \stackrel{d}{=} \BmH (t_k).
\end{eqnarray*}

\emph{(b)} $\BmH$ is approximately \emph{self-similar} in the following
sense. If
$a = 2^{2m_0}$, where $m_0$ is an integer, $m_0 \ge -m$, then for any
$k$ non-negative integer for which $ka$ is also an integer one has that
\begin{eqnarray*}
\lefteqn{a^{-H}\BmH (a k 2^{-2m})} \\
&=& \frac{a^{-H} 2^{-2Hm}}{\Gaa}
\sum_{r=-\infty}^{ak-1} \lsb (ak-r)^{\Hh} - (-r)_+^{\Hh} \rsb
\tilde X_m(r+1) \\
&=& \frac{2^{-2H(m+m_0)}}{\Gaa} \sum_{r=-\infty}^{k 2^{2m_0}-1} \lsb
(k 2^{2m_0}-r)^{\Hh} - (-r)_+^{\Hh} \rsb \tilde X_m(r+1) \\
&\stackrel{d}{=}& B_{m+m_0}^{(H)}(k 2^{-2m}). \label{eq:ssd1}
\end{eqnarray*}
On the other hand, Lemma~\ref{le:main} (and Theorem~\ref{th:main})
below show that $\BmH$ and $\BmpH$ (and $B_{m+n}^{(H)}$)
are uniformly close with arbitrary large probability on any
compact interval if $m$ is large enough (when $H > \quart$). It could
be proved in a similar fashion that for $a=j$, where $j \ge 0$ is an
arbitrary integer, $2^{2n} \le j \le 2^{2(n+1)}$ with an integer
$n \ge 0$, the finite dimensional distributions of
\[ a^{-H}\BmH (a k 2^{-2m}) = \frac{2^{-H(2m+\log_2 j)}}{\Gaa}
\sum_{r=-\infty}^{j k-1} \lsb (j k-r)^{\Hh} - (-r)_+^{\Hh} \rsb
\tilde X_m(r+1) \]
can be made arbitrarily close to the the finite dimensional
distributions of $B_{m+n}^{(H)}$ if $m$ is large
enough. Consequently, $\BmH$ is arbitrarily close to self-similar for
any dyadic $a=j2^{2m_0}$ if $m$ is large enough.

\emph{(c)} For any $0 < t_1 < \cdots < t_n$, the limit distribution of
the vector
\[ \lpa \BmH (t_1^{(m)}), \BmH (t_2^{(m)}), \dots , \BmH (t_n^{(m)})
\rpa \]
as $m \to \infty $ is \emph{Gaussian}, where  $t_j^{(m)}
= \lfloor t_j 2^{2m} \rfloor 2^{-2m}$, $1 \le j \le n$. This fact
follows from Theorem~\ref{th:main} (based on Lemma~\ref{le:Gauss})
below that states that the process $\BmH$ almost surely converges
to the Gaussian process $\WH$ on compact intervals.


\section{Convergence of the approximation to fBM}

At first it will be shown that two consecutive approximations of fBM
defined by (\ref{eq:afbrow1}), or equivalently by (\ref{eq:afbrow2}),
are uniformly close if $m$ is large enough, supposing $H > \quart$.
Apparently, the above RW approximation of BM is not good enough
to have convergence for $H \le \quart$.

When proving convergence, a large deviation inequality similar to
Lemma~\ref{le:basic} will play an important role.  If $X_1, X_2,\dots $
is a sequence of i.i.d. random variables, $\Pb X_k = \pm 1 \rb = \half$,
and $S = \sum_r a_r X_r$, where not all $a_r \in \mathbb{R}$ are zero and
$\Var (S) = \sum_r a_r^2 < \infty$, then
\begin{equation}
\Pb |S| \ge x \lpa \Var (S)\rpa ^{\half} \rb \le 2
e^{-\frac{x^2}{2}} \qquad (x \ge 0), \label{eq:ldp}
\end{equation}
see e.g. p. 33 in \cite{Str93}. The summation above may extend either
to finitely many or to countably many terms.

As a corollary, if $S_1, S_2, \dots , S_N$ are arbitrary sums of the
above type, one can get the following analog of Lemma~\ref{le:basic}.
For any $C >1$ and $N \ge 1$,
\begin{eqnarray}
\lefteqn{\Pb \max_{1\le k \le N} |S_k|
\ge (2C\log N)^{\half} \max_{1\le k \le N} \lpa \Var (S_k)\rpa ^\half
\rb} \nonumber \\
&\le& \sum_{k=1}^N \Pb |S_k| \ge (2C\log N \: \Var (S_k))^\half
\rb \le 2 N e^{-C \log N} = 2 N^{1-C}. \label{eq:basic2}
\end{eqnarray}

\begin{lem} \label{le:main}
For any $H \in (\quart, 1)$, $K > 0$, $C \ge 3$, and $m \ge m_3(C)$,
we have
\[ \Pb \max_{0 \le t_k \le K} |\BmpH (t_k) -
 \BmH (t_k)| \ge \alpha (H,K) m 2^{-\beta (H) m} \rb
\le 8 (K 2^{2m})^{1-C}, \]
where $t_k = k 2^{-2m}$ for $k \ge 0$ integers, $\beta(H) = \min(2H
- \half , \half )$ and
\[\alpha (H,K)
= \frac{(\log_* K)^{\half}}{\Gaa} \lsb \frac{\la \Hh \ra}{(1-H)^{\half}}
+ (\log_* K)^{\quart} \lpa 8 K^{\quart} + 36 \la \Hh \ra K^{H-\quart}
\rpa \rsb \]
if $H \in \lpa \quart , \half \rpa$,
\[\alpha (H,K)
= \frac{(\log_* K)^{\half}}{\Gaa} \lsb \frac{\la \Hh \ra}{(1-H)^{\half}}
+ (\log_* K)^{\quart} \lpa 5 + 312 \la \Hh \ra \rpa
K^{H-\quart} \rsb \]
if $H \in \lpa \half , 1\rpa$.
(The case $H=\half$ is described by Lemma~\ref{le:apprb}.)
\end{lem}

\begin{proof}
The proof is long, but elementary. Introduce the
following abbreviations: $\Delta B_m(t) = B_m(t+\Delta t)
- B_m(t)$, $\Delta B_{m+1}(t) = B_{m+1}(t+\quart\Delta t) - B_{m+1}(t)$.
Using (\ref{eq:afbrow1}) and then substituting $u=4r+j$, one gets that
\begin{eqnarray*}
\lefteqn{\BmpH (t_k) = \BmpH (4k 2^{-2(m+1)})} \\
&=& \frac{2^{-2H(m+1)}}{\Gaa} \sum_{u=-\infty}^{4k-1} \lsb (4k-u)^{\Hh}
- (-u)_+^{\Hh} \rsb \tilde X_{m+1}(u+1) \\
&=& \frac{2^{-2Hm-1}}{\Gaa} \sum_{r=-\infty}^{k-1} \sum_{j=0}^3 \lsb
\lpa k-r-\frac{j}{4}\rpa ^{\Hh}
- \lpa -r-\frac{j}{4}\rpa _+^{\Hh} \rsb  \\
&& \times \tilde X_{m+1}(4r+j+1) \\
&=&  \Ga \sum_{r=-\infty}^{k-1} \sum_{j=0}^3 \lsb (t_k-\trj )
^{\Hh} - (-\trj )_+^{\Hh} \rsb \Delta B_{m+1}(\trj ).
\end{eqnarray*}
So, subtracting and adding a suitable ``intermediate'' term, one
arrives at
\begin{eqnarray}
\lefteqn{\Gaa \lsb \BmpH (t_k) - \BmH (t_k) \rsb} \nonumber \\
&=& \sum_{r=-\infty}^{k-1} \sum_{j=0}^3 \lsb
(t_k-\trj )^{\Hh} - (-\trj )_+^{\Hh} \rsb
\Delta B_{m+1}(\trj ) \nonumber \\
&& - \lsb (t_k-t_r)^{\Hh} - (-t_r)_+^{\Hh} \rsb
\quart \Delta B_m(t_r) \nonumber \\
&=& \sum_{r=-\infty}^{k-1} \sum_{j=0}^3 \lb \lsb
(t_k-\trj )^{\Hh} - (-\trj )_+^{\Hh} \rsb \right. \nonumber \\
&& \left. - \lsb (t_k-t_r)^{\Hh} - (-t_r)_+^{\Hh} \rsb \rb
\Delta B_{m+1}(\trj ) \nonumber \\
&& + \sum_{r=-\infty}^{k-1} \sum_{j=0}^3 \lsb (t_k-t_r)^{\Hh}
- (-t_r)_+^{\Hh} \rsb  \lsb \Delta B_{m+1}(\trj )
- \quart \Delta B_m(t_r) \rsb \nonumber \\
&=:& (Z_{m,k} + Y_{m,k} + V_{m,k} + U_{m,k}). \label{eq:ZYVU1}
\end{eqnarray}

Here we introduced the following notations:
\begin{multline}
Z_{m,k} = \sum_{r=0}^{k-1} \sum_{j=0}^3
\lsb (t_k-\trj )^{\Hh} - (t_k-t_r)^{\Hh} \rsb
\Delta B_{m+1}(\trj ) \\
= 2^{-2Hm-1} \sum_{r=0}^{k-1} \sum_{j=0}^3 \lsb
(k-r-\frac{j}{4})^{\Hh} - (k-r)^{\Hh} \rsb \tilde X_{m+1}(4r+j+1)
\label{eq:Z} 
\end{multline}
and
\begin{eqnarray}
Y_{m,k} &= &\sum_{r=0}^{k-1} (t_k-t_r)^{\Hh}
\sum_{j=0}^3  \lsb \Delta B_{m+1}(\trj )
- \quart \Delta B_m(t_r) \rsb \nonumber \\
&=& \sum_{r=0}^{k-1} (t_k-t_r)^{\Hh}
\lb \lsb B_{m+1}(t_{r+1}) - B_{m+1}(t_r) \rsb - \lsb B_m(t_{r+1})
- B_m(t_r) \rsb \rb \nonumber \\
&=& \sum_{r=0}^{k} \lsb (t_k-t_{r-1})^{\Hh} - (t_k-t_r)^{\Hh} \rsb
\lsb B_{m+1}(t_r) - B_m(t_r) \rsb ,  \label{eq:Y}
\end{eqnarray}
applying ``summation by parts'' in the last row, as
in (\ref{eq:afbrow2}). Similarly, we introduced the following notations
for the corresponding ``infinite parts'' in (\ref{eq:ZYVU1}) (using
$v=-r$):
\begin{multline}
V_{m,k} \\
= 2^{-2Hm-1} \sum_{v=1}^{\infty} \sum_{j=0}^3
\lsb (k+v-\frac{j}{4})^{\Hh} - (v-\frac{j}{4})^{\Hh}
 - (k+v)^{\Hh} + v^{\Hh} \rsb  \\
\times \tilde X_{m+1}(-4v+j+1), \label{eq:V}
\end{multline}
and
\begin{eqnarray}
U_{m,k} &=&\sum_{v=1}^{\infty} \lsb (t_k+t_v)^{\Hh} - t_v^{\Hh} \rsb
\sum_{j=0}^3  \lsb \Delta B_{m+1}(-t_{v-\frac{j}{4}})
- \quart \Delta B_m(-t_v) \rsb \nonumber \\
&=& \sum_{v=1}^{\infty} \lsb (t_k+t_v)^{\Hh} - t_v^{\Hh} \rsb
\nonumber \\
&& \times \lb \lsb B_{m+1}(-t_{v-1}) - B_{m+1}(-t_v) \rsb
- \lsb B_m(-t_{v-1}) - B_m(-t_v) \rsb \rb \nonumber \\
&=& \sum_{v=1}^{\infty} \lsb (t_k+t_{v+1})^{\Hh} - t_{v+1}^{\Hh}
- (t_k+t_v)^{\Hh} + t_v^{\Hh} \rsb \nonumber \\
&& \times \lsb B_{m+1}(-t_v) - B_m(-t_v) \rsb . \label{eq:U}
\end{eqnarray}

The maxima of $Z_{m,k}, Y_{m,k}, V_{m,k}$ and $U_{m,k}$ can be estimated
separately:
\begin{eqnarray}
\lefteqn{\max_{0 \le t_k \le K} |\BmpH (t_k) - \BmH (t_k)|} \\
&\le& \Ga \lpa \max_k |Z_{m,k}| + \max_k |Y_{m,k}| + \max_k |V_{m,k}|
+ \max_k |U_{m,k}| \rpa ,  \nonumber \label{eq:ZYVU}
\end{eqnarray}
where each maximum on the right hand side is taken for
$1 \le k \le K 2^{2m}$ and one can suppose that $K 2^{2m} \ge 1$,
that is, $\Delta t \le K$, since otherwise the maximal difference
in (\ref{eq:ZYVU}) is zero.

\emph{(a) The maximum of $Z_{m,k}$.}
In the present case the large deviation inequality (\ref{eq:ldp}), or
rather, its corollary (\ref{eq:basic2}) is applied. By (\ref{eq:Z}),
\begin{eqnarray*}
\Var (Z_{m,k}) &=& 2^{-4Hm-2} \sum_{r=0}^{k-1} \sum_{j=0}^3
\lsb (k-r-\frac{j}{4})^{\Hh} - (k-r)^{\Hh} \rsb ^2 \\
&=& 2^{-4Hm-2} \sum_{r=0}^{k-1} \sum_{j=0}^3 (k-r)^{2H-1} \lsb
(1 - \frac{j}{4(k-r)})^{\Hh} - 1 \rsb ^2 .
\end{eqnarray*}

The term in brackets can be estimated using a binomial series with
$0 \le j \le 3$, $k-r \ge 1$:
\begin{eqnarray*}
&&\la \lpa 1 - \frac{j}{4(k-r)} \rpa ^{\Hh} - 1 \ra = \la
\sum_{s=1}^{\infty} \binom{\Hh}{s} (-1)^s \lpa \frac{j}{4(k-r)}
\rpa ^s \ra \\
&\le& \la \Hh \ra \frac{j}{4(k-r)} \lpa 1-\frac{j}{4(k-r)} \rpa ^{-1}
\le \la \Hh \ra \frac{j}{4(k-r)} \lpa 1-\frac{j}{4} \rpa ^{-1}.
\end{eqnarray*}
Thus
\[ \sum_{j=0}^3 \lsb \lpa 1 - \frac{j}{4(k-r)}\rpa ^{\Hh} - 1 \rsb ^2
\le \lpa \Hh \rpa ^2 ~\frac{91}{9}~ \frac{1}{(k-r)^2}. \]
Also, if $k \ge 1$, one has
\[ \sum_{r=0}^{k-1} (k-r)^{2H-3} < 1 + \int_0^{k-1} (k-x)^{2H-3} \di x
\le \threehalf ~\frac{1}{1-H} . \]
Then for any $k \ge 1$ it follows that
\[ \Var (Z_{m,k}) \le 2^{-4Hm} \lpa \Hh \rpa ^2 \frac{273}{72}~
\frac{1}{1-H}. \]
Hence taking $N=K 2^{2m}$ and $C > 1$ in (\ref{eq:basic2}), one obtains
that
\begin{eqnarray*}
\lefteqn{\Pb \max_{1 \le k \le N} |Z_{m,k}| \ge \lpa \frac{273}{72} \rpa
^{\half}~ \la \Hh \ra (1-H)^{-\half} 2^{-2Hm} (2C\log N)^\half \rb}  \\
&\le& \sum_{k=1}^N \Pb \la Z_{m,k} \ra
\ge \lpa 2C\log N \: \Var (Z_{m,k}) \rpa ^\half \rb \le 2 N^{1-C}.
\end{eqnarray*}
Since
\begin{equation}
\log N = \log \lpa K 2^{2m} \rpa \le (1+\log 4) m \log_* K
\le 2.5 m \log_* K , \label{eq:logN}
\end{equation}
one obtains the following result:
\begin{equation}
\max_{1 \le k \le K 2^{2m}} |Z_{m,k}| \le 5 \la \Hh \ra (1-H)
^{-\half} (C \log_* K)^\half~m^{\half} 2^{-2Hm}, \label{eq:Zmax}
\end{equation}
with the exception of a set of probability at most $2 \lpa K 2^{2m} \rpa
^{1-C}$, where $m \ge 1$, $K > 0$ and $C > 1$ are arbitrary.

\emph{(b) The maximum of $Y_{m,k}$.}
By its definition (\ref{eq:Y}),
\begin{eqnarray*}
\max_{1 \le k \le K 2^{2m}} |Y_{m,k}| &\le& \max_{0 \le t_r \le K}
\la B_{m+1}(t_r) - B_m(t_r) \ra \\
&& \times \max_{0 \le t_k \le K} \sum_{0 \le t_r
\le t_k} \la (t_k-t_{r-1})^{\Hh} - (t_k-t_r)^{\Hh} \ra .
\end{eqnarray*}
The first factor, the maximal difference between two consecutive
approximations of BM appearing here can be estimated by
Lemma~\ref{le:apprb}. For the second factor one can apply a binomial
series:
\begin{eqnarray*}
\lefteqn{\sum_{r=0}^k \la (t_k-t_{r-1})^{\Hh} - (t_k-t_r)^{\Hh} \ra} \\
&=& 2^{-m(2H-1)} \lb 1 + \la 2^{\Hh} - 1 \ra + \sum_{r=0}^{k-2}
(k-r)^{\Hh} \la \lpa 1 + \frac{1}{k-r} \rpa ^{\Hh} - 1 \ra \rb \\
&=& 2^{-m(2H-1)} \lb 1 + \la 2^{\Hh} - 1 \ra + \sum_{r=0}^{k-2}
(k-r)^{\Hh} \la \sum_{s=1}^{\infty} \binom{\Hh}{s} ~ \frac{1}{(k-r)^s}
\ra \rb \\
&\le& 2^{-m(2H-1)} \lb 1 + \la 2^{\Hh} - 1 \ra + \sum_{r=0}^{k-2}
(k-r)^{\Hh} \frac{\la \Hh \ra}{k-r} \rb .
\end{eqnarray*}

Since for $H \ne \half$
\[ \sum_{r=0}^{k-2} (k-r)^{\Hth} \le \int_0^{k-1}
(k-x)^{\Hth} \di x = \frac{1-k^{\Hh}}{\half - H}, \]
it follows for any $m \ge 0$ that
\begin{eqnarray*}
\lefteqn{\max_{1 \le k \le K 2^{2m}} \sum_{r=0}^k \la
(t_k-t_{r-1})^{\Hh} - (t_k-t_r)^{\Hh} \ra} \\
&\le& 2^{-m(2H-1)} \max_{1 \le k \le K 2^{2m}} \lb 1 + \la 2^{\Hh} - 1
\ra + \la 1 - k^{\Hh} \ra \rb \\
&\le& \left\{
\begin{array}{ll}
2^{-2m(\Hh)} \lpa 3-2^{\Hh} \rpa \le 3~ 2^{-2m(\Hh)} & \mbox{ if }
0 < H < \half , \\
2^{-2m(\Hh)} \lpa 2^{\Hh} - 1 \rpa + K^{\Hh} \le (2K)^{\Hh} &
\mbox{ if } \half < H < 1.
\end{array} \right.
\end{eqnarray*}
(In the last row we used that here $2^{-2m} \le K$.)

Combining this with Lemma~\ref{le:apprb}, we obtain the result
\begin{equation}
\max_{1 \le k \le K 2^{2m}} |Y_{m,k}| \le \left\{
\begin{array}{ll}
3 K^{\quart} (\log_*K)^{\frac{3}{4}} m 2^{-2m\lpa H - \quart \rpa}&
\mbox{ if } \quart < H < \half , \\
2^{\Hh} K^{H-\quart} (\log_*K)^{\frac{3}{4}} m 2^{-\frac{m}{2}}&
\mbox{ if } \half < H < 1.
\end{array} \right. \label{eq:Ymax}
\end{equation}
with the exception of a set of probability at most $3 \lpa K 2^{2m} \rpa
^{1-C}$, where $K > 0$, $C > 1$ are arbitrary, and $m \ge m_1(C)$.
Thus in the case $0<H<\half$ we have only a partial result:
the relative weakness of the above-described RW approximation of BM
causes that apparently we have no convergence for $0<H\le \quart$.

\emph{(c) The maximum of $V_{m,k}$.}
Here one can use the same idea as in part (a), including
the application of the corollary (\ref{eq:basic2}) of the large
deviation principle. We begin with (\ref{eq:V}),
\begin{eqnarray*}
\lefteqn{\Var (V_{m,k})} \\
&=& 2^{-4Hm-2} \sum_{v=1}^{\infty} \sum_{j=0}^3
\lsb (k+v-\frac{j}{4})^{\Hh} - (k+v)^{\Hh} - (v-\frac{j}{4})^{\Hh}
+ v^{\Hh} \rsb ^2 \\
&=& 2^{-4Hm-2} \sum_{v=1}^{\infty} \sum_{j=0}^3 \lb
(k+v)^{\Hh} \lsb \lpa 1-\frac{j}{4(k+v)} \rpa ^{\Hh} - 1 \rsb \right. \\
&& \left. - v^{\Hh} \lsb \lpa 1-\frac{j}{4v} \rpa ^{\Hh} - 1 \rsb \rb ^2 .
\end{eqnarray*}
As in (a), now we use binomial series for the expressions in brackets
($k \ge 1$, $0 \le j \le 3$, $v \ge 1$):
\begin{eqnarray*}
A &=& (k+v)^{\Hh} \lsb \lpa 1-\frac{j}{4(k+v)} \rpa ^{\Hh} - 1 \rsb \\
&=& (k+v)^{\Hh} \sum_{s=1}^{\infty} \binom{\Hh}{s} (-1)^s \lpa
\frac{j}{4(k+v)} \rpa ^s
\end{eqnarray*}
and
\[ B = v^{\Hh} \lsb \lpa 1-\frac{j}{4v} \rpa ^{\Hh} - 1 \rsb =
v^{\Hh} \sum_{s=1}^{\infty} \binom{\Hh}{s} (-1)^s \lpa \frac{j}{4v}
\rpa ^s. \]
Then $A$ and $B$ have the same sign and $0 \le A_1 \le |A| \le |B| \le
B_2$, where
\[ A_1 = (k+v)^{\Hh} \la \Hh \ra \frac{j}{4(k+v)}
= \la \Hh \ra \frac{j}{4} (k+v)^{\Hth} \]
and
\[ B_2 = v^{\Hh} \la \Hh \ra \frac{j}{4v} \lpa 1 - \frac{j}{4v}
\rpa ^{-1} \le \la \Hh \ra \frac{j}{4-j} v^{\Hth}. \]
Hence
\[ (A-B)^2 \le A_1^2 + B_2^2 \le \lpa \Hh \rpa ^2 \lsb \lpa \frac{j}{4}
\rpa ^2 (k+v)^{2H-3} + \lpa \frac{j}{4-j} \rpa ^2 v^{2H-3} \rsb . \]
Since for any $k \ge 0$,
\[ \sum_{v=1}^{\infty} (k+v)^{2H-3} < 1 + \int_1^{\infty} (k+x)^{2H-3}
\di x = 1 + \frac{(k+1)^{2H-2}}{2-2H} < \threehalf \frac{1}{1-H} \]
it follows that
\begin{multline*}
\Var (V_{m,k}) \\
\le  2^{-4Hm-2} \lpa \Hh \rpa ^2
\sum_{v=1}^{\infty} \sum_{j=0}^3 \lsb \lpa \frac{j}{4}
\rpa ^2 (k+v)^{2H-3}
+ \lpa \frac{j}{4-j} \rpa ^2 v^{2H-3} \rsb \\
\le \frac{791}{192} ~\frac{\lpa \Hh \rpa ^2}{1-H} 2^{-4Hm} .
\end{multline*}
Applying corollary (\ref{eq:basic2}) of the large deviation inequality
with $N=K 2^{2m}$ one obtains that
\begin{multline*}
\Pb \max_{1 \le k \le N} |V_{m,k}| \ge \lpa \frac{791}{192}
\rpa ^\half~ \la \Hh \ra (1-H)^{-\half} 2^{-2Hm} (2C\log N)^\half \rb
\\
\le \sum_{k=1}^N \Pb \la V_{m,k} \ra
\ge \lpa 2C\log N \: \Var (V_{m,k}) \rpa ^\half \rb \le 2 N^{1-C}.
\end{multline*}
Hence using (\ref{eq:logN}) one gets the result
\begin{equation}
\max_{1 \le k \le K 2^{2m}} |V_{m,k}| \le 5 \la \Hh \ra (1-H)
^{-\half} (C \log_* K)^\half~m^{\half} 2^{-2Hm},
\label{eq:Vmax}
\end{equation}
with the exception of a set of probability at most $2 \lpa K 2^{2m} \rpa
^{1-C}$, where $m \ge 1$, $K > 0$ and $C > 1$ are arbitrary.

\emph{(d) The maximum of $U_{m,k}$.}
We divide the half line into intervals of length $L$, where $L \ge
4K$. For definiteness, choose $L=4K$. Apart from this,
this part will be similar to part (b). In the sequel we use the
convention that when the lower limit of a summation is a real number
$x$, the summation starts at $\lceil x \rceil$, and similarly, if the
upper limit is $y$, the summation ends at $\lfloor y \rfloor$.
By (\ref{eq:U}),
\begin{multline}
|U_{m,k}|
\le \sum_{j=1}^{\infty} \sum_{(j-1)L < t_v \le jL}
\la (t_k+t_{v+1})^{\Hh} - (t_k+t_v)^{\Hh} - t_{v+1}^{\Hh}
+ t_v^{\Hh} \ra \\
\qquad \times \la B_{m+1}(-t_v) - B_m(-t_v) \ra  \\
\le \sum_{j=1}^{\infty} \max_{(j-1)L < t_v \le jL}
\la B_{m+1}(-t_v) - B_m(-t_v) \ra  (\Delta t)^{\Hh} \label{eq:Ufact} \\
\times  \sum_{v=(j-1)L2^{2m}+1}^{jL2^{2m}}
\la (k+v+1)^{\Hh} - (k+v)^{\Hh} - (v+1)^{\Hh} + v^{\Hh} \ra  .
\end{multline}

Lemma~\ref{le:apprb} gives an upper bound for the maximal
difference between two consecutive approximations of BM if $j \ge
1$ is an arbitrary fixed value:
\begin{eqnarray}
\lefteqn{\max_{(j-1)L<t_v\le jL} \la B_{m+1}(-t_v) - B_m(-t_v) \ra} \\
&\le& (jL)^{\quart} \lpa \log_*(jL) \rpa ^{\frac{3}{4}} m 2^
{-\frac{m}{2}} \nonumber \\
&\le& \left\{
\begin{array}{ll}
L^{\quart} (\log_*L)^{\frac{3}{4}} m 2^{-\frac{m}{2}} & \mbox{ if } j=1,
 \\
2 j^{\quart}(\log_*j)^{\frac{3}{4}} L^{\quart} (\log_*L)^{\frac{3}{4}}
m 2^{-\frac{m}{2}} & \mbox{ if } j \ge 2,
\end{array} \right. \label{eq:Ufact1}
\end{eqnarray}
with the exception of a set of probability at most $3 \lpa jL 2^{2m}
\rpa ^{1-C}$, where $C > 1$ is arbitrary and $m \ge m_1(C)$. This
implies for any $C \ge 3$ and $m \ge m_1(C)$ that the above inequality
(\ref{eq:Ufact1}) holds simultaneously for all $j=1,2,3, \dots$ with the
exception of a set of probability at most
\begin{equation}
3 \lpa L 2^{2m} \rpa ^{1-C} \sum_{j=1}^{\infty} j^{1-C} < 3 \lpa L
2^{2m} \rpa ^{1-C} \frac{\pi ^2}{6} < \lpa K 2^{2m} \rpa ^{1-C}.
\label{eq:Uprob}
\end{equation}

For the other major factor in (\ref{eq:Ufact}) binomial series are
applied as above, with $m \ge 0$, $k \ge 1$, and $v \ge 1$:
\begin{eqnarray*}
A &=& (k+v+1)^{\Hh} - (k+v)^{\Hh} = (k+v)^{\Hh} \lsb \lpa 1+\frac{1}{k+v}
\rpa ^{\Hh} - 1 \rsb \\
&=& (k+v)^{\Hh} \sum_{s=1}^{\infty} \binom{\Hh}{s} \frac{1}{(k+v)^s} ,
\end{eqnarray*}
and for $v \ge 2$:
\[ B = (v+1)^{\Hh} - v^{\Hh}  =  v^{\Hh} \lsb \lpa 1+\frac{1}{v}
\rpa^{\Hh} - 1 \rsb = v^{\Hh} \sum_{s=1}^{\infty} \binom{\Hh}{s}
\frac{1}{v^s}, \]
while $B=2^{\Hh}-1$ when $v=1$.
Then $A$ and $B$ have the same sign, $0 \le A_1 \le |A| \le |B| \le B_2$,
and so $|A-B| \le B_2-A_1$, where
\[ A_1 = (k+v)^{\Hth} \la \Hh \ra - (k+v)^{H-\frac{5}{2}} ~\half ~
\la \Hh \ra \lpa \threehalf - H \rpa  \]
and $B_2 =v^{\Hth} \la \Hh \ra$ .

Thus if the second major factor in (\ref{eq:Ufact}) is denoted by
$C_{m,k,j}$, we obtain for any $j \ge 1$ that
\begin{multline*}
C_{m,k,j} \\
= (\Delta t)^{\Hh} \sum_{v=(j-1)L2^{2m}+1}^{jL2^{2m}}
\la (k+v+1)^{\Hh} - (k+v)^{\Hh} - (v+1)^{\Hh} + v^{\Hh} \ra \\
\le \la \Hh \ra (\Delta t)^{\Hh} \\
\times \sum_{v=(j-1)L2^{2m}+1}^{jL2^{2m}}\lb v^{\Hth}
- (k+v)^{\Hth} + \half\, \lpa \threehalf - H \rpa
(k+v)^{H-\frac{5}{2}} \rb .
\end{multline*}
For $H \ne \half$ one can get the estimates for $j=1$:
\[ \sum_{v=1}^{L2^{2m}} v^{\Hth} < 1 + \int_1
^{L2^{2m}} x^{\Hth} \di x = \frac{(\Delta t)^{\half -H}}{\Hh} L^{\Hh}
+ \frac{\Hth}{\Hh}, \]
and for $j \ge 2$:
\begin{multline*}
\sum_{v=(j-1)L2^{2m}+1}^{jL2^{2m}} v^{\Hth} 
< \int_{(j-1)L2^{2m}}^{jL2^{2m}} x^{\Hth} \di x \\
= \frac{(\Delta t)^{\half -H}}{\Hh} \lsb
(jL)^{\Hh} - \lpa (j-1)L \rpa ^{\Hh} \rsb,
\end{multline*}
further, for any $j \ge 1$,
\begin{eqnarray*}
\lefteqn{\sum_{v=(j-1)L2^{2m}+1}^{jL2^{2m}}(k+v)^{\Hth}
> \int_{(j-1)L2^{2m}+1}^{jL2^{2m}+1} (k+x)^{\Hth} \di x} \\
&=& \frac{(\Delta t)^{\half -H}}{\Hh} \lsb (t_{k+1}+jL)
^{\Hh} - \lpa t_{k+1}+(j-1)L \rpa ^{\Hh} \rsb,
\end{eqnarray*}
and also for any $j \ge 1$,
\begin{multline*}
\sum_{v=(j-1)L2^{2m}+1}^{jL2^{2m}}(k+v)^{H-\frac{5}{2}}
< \int_{(j-1)L2^{2m}}^{jL2^{2m}} (k+x)^{H-\frac{5}{2}} \di x \\
= \frac{(\Delta t)^{\threehalf-H}}{\threehalf-H} \lsb \lpa
t_k+(j-1)L \rpa ^{\Hth} - (t_k+jL)^{\Hth} \rsb.
\end{multline*}
Denote the sign of a real number $x$ by $\epsilon_x$ ($0$ if $x=0$).
When $j=1$, it follows that
\begin{eqnarray*}
C_{m,k,1} &\le&
 \epsilon_{\Hh} \lb L^{\Hh} \lsb 1 - \lpa 1+\frac{t_{k+1}}{L}\rpa
^{\Hh} + \lpa 1+\frac{t_{k+1}-L}{L}\rpa ^{\Hh} \rsb \right. \\
&& \left. + \lpa \Hth \rpa \lpa \Delta t \rpa ^{\Hh} \rb \\
&& + \la \Hh \ra \frac{\Delta t}{2}
L^{\Hth} \lsb \lpa 1+\frac{t_k-L}{L}\rpa ^{\Hth} - \lpa
1+\frac{t_k}{L}\rpa ^{\Hth} \rsb ,
\end{eqnarray*}
and similarly, when $j \ge 2$,
\begin{eqnarray*}
\lefteqn{C_{m,k,j} \le \epsilon_{\Hh} \lsb (jL)^{\Hh} - \lpa (j-1)L
\rpa^{\Hh} \right.} \\
&& \left. - (t_{k+1}+jL)^{\Hh} + \lpa t_{k+1}+(j-1)L \rpa ^{\Hh}
\rsb \\
&& + \la \Hh \ra \frac{\Delta t}{2} \lsb \lpa t_k+(j-1)L \rpa ^{\Hth}
- (t_k+jL)^{\Hth} \rsb \\
&=& \epsilon_{\Hh} (jL)^{\Hh} \\
&&\times \lsb 1-\lpa 1-\frac{1}{j}\rpa
^{\Hh} - \lpa 1+\frac{t_{k+1}}{jL}\rpa ^{\Hh} + \lpa
1+\frac{t_{k+1}-L}{jL}\rpa ^{\Hh} \rsb \\
&& + \la \Hh \ra \frac{\Delta t}{2}
(jL)^{\Hth} \lsb \lpa 1+\frac{t_k-L}{jL}\rpa ^{\Hth} - \lpa
1+\frac{t_k}{jL}\rpa ^{\Hth} \rsb ,
\end{eqnarray*}

Applying binomial series here again, first we get when $j \ge 2$ that
\begin{eqnarray*}
\lefteqn{\epsilon_{\Hh} \lsb 1-\lpa 1-\frac{1}{j}\rpa^{\Hh} \rsb
= \epsilon_{\Hh} \sum_{s=1}^{\infty} \binom{\Hh}{s} \frac{(-1)
^{s+1}}{j^s}} \\
&\le& \la \Hh \ra \frac{1}{j} + \half ~ \la \Hh \ra \lpa \threehalf-H
\rpa \frac {1}{j^2} \lpa 1-\frac{1}{j} \rpa ^{-1} \\
&\le& \la \Hh \ra \frac{1}{j} + \la \Hh \ra \lpa \threehalf-H
\rpa \frac {1}{j^2},
\end{eqnarray*}
since each term of the series is positive. Furthermore, with any
$j \ge 1$,
\begin{eqnarray*}
\lefteqn{\epsilon_{\Hh} \lsb \lpa 1-\frac{L-t_{k+1}}{jL}\rpa
^{\Hh} - \lpa 1+\frac{t_{k+1}}{jL}\rpa ^{\Hh} \rsb} \\
&=& \epsilon_{\Hh} \sum_{s=1}^{\infty} \binom{\Hh}{s}
\frac{(-1)^s} {(jL)^s} \lsb (L-t_{k+1})^s - (-t_{k+1})^s \rsb \le -\la
\Hh \ra \frac{1}{j},
\end{eqnarray*}
since each term of the series is negative: $L=4K
\ge 2 t_{k+1}$, and the term in brackets is not larger than
$2(L-\Delta t)^s$. Finally,
\begin{multline*}
\lpa 1-\frac{L-t_k}{jL}\rpa
^{\Hth} - \lpa 1+\frac{t_k}{jL}\rpa ^{\Hth} \\
= \sum_{s=1}^{\infty} \binom{\Hth}{s} \frac{(-1)^s} {(jL)^s} \lsb
(L-t_k)^s - (-t_k)^s \rsb \\
\le \lpa \threehalf-H \rpa \frac{4}{3j} \lpa 1-\frac{L-\Delta t}{jL}
\rpa ^{-1} \le \lpa \threehalf-H \rpa \frac{4L}{3j\Delta t},
\end{multline*}
since each term of the series is positive and the term in brackets is
not larger than $\frac43(L-\Delta t)^s$.

Thus when $j \ge 2$ it follows for any $m \ge 0$, $k \ge 1$ that
\begin{eqnarray*}
\lefteqn{C_{m,k,j}} \\
&\le& (jL)^{\Hh} \la \Hh \ra \lpa \threehalf-H \rpa
\frac {1}{j^2}
+ \la \Hh \ra \frac{\Delta t}{2} (jL)^{\Hth} \lpa \threehalf-H \rpa
\frac{4L}{3j\Delta t} \\
&\le& a_H \la \Hh \ra L^{\Hh} j^{H-\frac52} \qquad \mbox{ where }
a_H = \left\{
\begin{array}{ll}
\frac52 & \mbox{ if } 0 < H < \half , \\
\frac53 & \mbox{ if } \half < H < 1 .
\end{array} \right.
\end{eqnarray*}
In a similar manner, when $j=1$ one can get for any $m \ge 0$, $k \ge 1$
that
\begin{eqnarray*}
C_{m,k,1} &\le& \epsilon_{\Hh} L^{\Hh} - \la \Hh \ra L^{\Hh}
+ \epsilon_{\Hh} \lpa \Hth \rpa (\Delta t)^{\Hh} \\
&& +  \la \Hh \ra \frac{\Delta t}{2} L^{\Hth} \lpa \threehalf-H \rpa
\frac{4L}{3\Delta t} \\
&=& \epsilon_{\Hh} \lpa \threehalf-H \rpa \lsb \frac{2}{3} H
L^{\Hh} - (\Delta t)^{\Hh} \rsb  \\
&\le& \left\{
\begin{array}{ll}
\threehalf (\Delta t)^{\Hh} & \mbox{ if } 0 < H < \half , \\
\frac38 L^{\Hh} & \mbox{ if } \half < H < 1 .
\end{array} \right.
\end{eqnarray*}

Then combine these results with (\ref{eq:Ufact1}) and (\ref{eq:Uprob})
in (\ref{eq:Ufact}). Using
\begin{eqnarray}
\lefteqn{\sum_{j=2}^{\infty} j^{\quart}(\log_*j)^{\frac{3}{4}}
j^{H -\frac{5}{2}} < \int_1^{\infty} x^{H-\frac{9}{4}} \log_* x \di x }\nonumber \\
&=& \int_1^e x^{H-\frac{9}{4}} \di x + \int_e^{\infty} x^{H-\frac{9}{4}}
\log x \di x  \nonumber \\
&=& \lpa \frac{5}{4}-H \rpa ^{-1} + \lpa \frac{5}{4}-H \rpa ^{-2}
e^{H-\frac{5}{4}} \nonumber \\
&<& \left\{
\begin{array}{rl}
2.5 & \mbox{ if } 0 < H < \half , \\
16.5 & \mbox{ if } \half < H < 1 ,
\end{array} \right. \label{eq:Uj}
\end{eqnarray}
one can get the result of part (d). Consider first the case $\half <H<1$:
\begin{eqnarray}
\lefteqn{\max_{1 \le k \le K 2^{2m}} |U_{m,k}|} \nonumber \\
&\le& L^{\quart} (\log_*L)^{\frac{3}{4}} m 2^{-\frac{m}{2}} \frac{3}{8}
L^{\Hh} + 33~ L^{\quart} (\log_*L)^{\frac{3}{4}} m 2^{-\frac{m}{2}}
\frac{5}{3}~ \la \Hh \ra L^{\Hh} \nonumber \\
&\le& \lpa 3+ 312 ~ \la \Hh \ra \rpa K^{H-\quart} (\log_*K)
^{\frac{3}{4}} m 2^{-\frac{m}{2}}, \label{eq:Umax1}
\end{eqnarray}
for any $C \ge 3$ and $m \ge m_1(C)$ with the exception of a set of
probability at most $\lpa K 2^{2m} \rpa ^{1-C}$. (Recall that $L=4K$.)

In the second case when $0<H<\half$ the above method apparently gives
convergence here (just like in part (b)) only when $\quart <H<\half$:
\begin{eqnarray}
\lefteqn{\max_{1 \le k \le K 2^{2m}} |U_{m,k}|} \nonumber \\
&\le& L^{\quart} (\log_*L)^{\frac34} m 2^{-\frac{m}{2}}
\threehalf(\Delta t)^{\Hh} + 5~ L^{\quart} (\log_*L)^{\frac34} m
2^{-\frac{m}{2}} \, \frac52~\la \Hh \ra L^{\Hh}
\nonumber \\
&\le& 5 \, K^{\quart} (\log_*K)^{\frac{3}{4}} m 2^{-2m(H-\quart )}
+ 36 \, \la \Hh \ra  K^{H-\quart} (\log_*K)^{\frac{3}{4}}
m 2^{-\frac{m}{2}}, \label{eq:Umax2}
\end{eqnarray}
for any $C \ge 3$ and $m \ge m_1(C)$ with the exception of a set of
probability at most $\lpa K 2^{2m} \rpa ^{1-C}$.

Now one can combine the results of parts (a), (b), (c), and (d),
see (\ref{eq:ZYVU}), (\ref{eq:Zmax}), (\ref{eq:Ymax}),
(\ref{eq:Vmax}), (\ref{eq:Umax1}), (\ref{eq:Umax2}), to obtain the
statement of the lemma. Remember that the rate of convergence in parts
(a) and (c) is faster than the one in parts (b) and (d). Particularly,
observe that there is a factor $m$ in (b) and (d) which has a
counterpart $m^{\half}$ in (a) and (c). Since in the statement of this
lemma we simply replaced the faster converging factors by the slower
converging ones, the constant multipliers in (a) and (c) can be ignored
if $m$ is large enough. \qed
\end{proof}

It is simple to extend formula (\ref{eq:afbrow2}) of the $m$th
approximation $\BmH$ of fBM to real arguments $t$ by linear
interpolation, just like in the case of the $m$th approximation $B_m(t)$
of ordinary BM, see e.g. in \cite{Sza96}. So let $m \ge 0$ and
$k\ge 0$ be integers, $\gamma \in [0, 1]$, and define that
\begin{eqnarray}
\lefteqn{\BmH (t_{k+\gamma}) = \gamma \BmH (t_{k+1}) + (1-\gamma )
\BmH (t_k)} \label{eq:intpol} \\
&=& \Ga \sum_{r=-\infty }^k \lsb (t_k - t_{r-1})^{\Hh}
- (t_k - t_r)^{\Hh} \rsb B_m(t_{r+\gamma}) \nonumber \\
&& + \lsb (-t_r)_+^{\Hh} - (- t_{r-1})_+^{\Hh} \rsb B_m(t_r).
\nonumber
\end{eqnarray}
Then the resulting \emph{continuous parameter approximations of fBM}
$\BmH (t)$ $(t \ge 0)$ have continuous, piecewise linear sample paths.
With this definition we are ready to state a main result of this paper.

\begin{thmn} \label{th:main}
For any $H \in (\quart, 1)$, the sequence $\BmH (t)$  $(t \ge 0, m=0,1,2,
\dots )$ a.s. uniformly converges to a fBM $\WH (t)$  $(t \ge 0)$ on any
compact interval $[0, K]$, $K > 0$. If $K > 0$, $C \ge 3$, and
$m \ge m_4(C)$, it follows that
\[ \Pb \max_{0 \le t \le K} |\WH (t) - \BmH (t)|
\ge \frac{\alpha (H,K)}{\lpa 1-2^{-\beta (H)} \rpa ^2} m 2^{-\beta (H) m}
\rb \le 9 (K 2^{2m})^{1-C}, \]
where  $\alpha (H,K)$ and $\beta (H) = \min(2H - \half , \half )$ are
the same as in Lemma~\ref{le:main}. (The case $H=\half$ is
described by Theorem~\ref{th:wiener}.)
\end{thmn}

\begin{proof}
At first we consider the maximum of $|\BmpH (t) - \BmH (t)|$ for real
$t \in [0, K]$. Lemma~\ref{le:main} gives an upper bound $D_m$ for their
maximal difference at vertices with $t=t_k=k \Delta t$:
\[ \max_{0 \le t_k \le K} |\BmpH (t_k) - \BmH (t_k)| \le D_m, \]
except for an event of probability at most $8(K2^{2m})^{1-C}$.
Since both $\BmpH (t)$ and $\BmH (t)$ have
piecewise linear sample paths, their maximal difference must occur at
vertices of the sample paths.  Let $M_m$ denote the maximal increase of
$\BmH$ between pairs of points $t_k, t_{k+1}$ in $[0, K]$:
\[ \max_{0 \le t_k \le K} |\BmH (t_{k+1}) - \BmH (t_k)| \le M_m, \]
except for an event of probability at most $2(K2^{2m})^{1-C}$,
cf. (\ref{eq:Mmax}) below. A sample path of $\BmpH (t)$ makes four steps
on any interval $[t_k, t_{k+1}]$. To compute its maximal deviation from
$D_m$ it is enough to estimate its change between the midpoint and an
endpoint of such an interval, at two steps from both the left and right
endpoints:
\[ \max_{0 \le t_k \le K} |\BmpH (t_{k \pm \half}) - \BmpH (t_k)|
\le 2 M_{m+1}, \]
except for an event of probability at most $2(K2^{2(m+1)})^{1-C}$.
Hence
\begin{eqnarray*}
\lefteqn{\max_{0 \le t_k \le K} |\BmpH (t_{k + \half})
- \BmH(t_{k + \half})|} \\
&=& \max_{0 \le t_k \le K} \la \BmpH (t_{k +\half}) - \textstyle\half
\lpa \BmH(t_k)+\BmH(t_{k+1})\rpa \ra \\
&\le& \max_{0 \le t_k \le K} |\BmpH (t_k) - \BmH (t_k)|
+ \max_{0 \le t_k \le K} |\BmpH (t_{k \pm \half}) - \BmpH(t_k)| \\
&\le& D_m + 2 M_{m+1},
\end{eqnarray*}
except for an event of probability at most
$(8+2^{3-2C})(K2^{2m})^{1-C}$. The explanation above
shows that at the same time this gives the upper bound we were looking
for:
\begin{equation}
\max_{0 \le t \le K} |\BmpH (t) - \BmH (t)| \le D_m + 2 M_{m+1},
\label{eq:BHdiff}
\end{equation}
except for an event of probability at most
$(8+2^{3-2C})(K2^{2m})^{1-C}$.

Thus we have to find an upper estimate $M_m$. For that the large
deviation inequality (\ref{eq:basic2}) will be used.
By (\ref{eq:afbrow1}), the increment of $\BmH (t)$ on $[t_k, t_{k+1}]$ is
\begin{eqnarray*}
A_{m,k} &=& |\BmH (t_{k+1}) - \BmH (t_k)| \\
&=& \frac{2^{-2Hm}}{\Gaa} \sum_{r=-\infty }^k \lsb (k+1-r)^{\Hh}
- (k-r)^{\Hh} \rsb \tilde X_m(r+1).
\end{eqnarray*}
Then a similar argument can be used as in the proof of
Lemma~\ref{le:main}, see e.g. part (a) there:
\begin{eqnarray*}
\lefteqn{\frac{\Gamma ^2 \lpa H+\half \rpa}{2^{-4Hm}} \Var (A_{m,k})
= \sum_{r=-\infty }^k \lsb (k+1-r)^{\Hh} - (k-r)^{\Hh} \rsb ^2} \\
&=& 1 + \lpa 2^{\Hh} - 1 \rpa ^2 + \sum_{r=-\infty }^{k-2} (k-r)^{2H-1}
\lsb \lpa 1 + \frac{1}{k-r} \rpa ^{\Hh} - 1 \rsb ^2 \\
&\le& 1 + \lpa 2^{\Hh} - 1 \rpa ^2 + \sum_{r=-\infty }^{k-2} (k-r)^{2H-1}
\lpa \Hh \rpa ^2 \frac{1}{(k-r)^2} \\
&\le& 1 + \lpa \Hh \rpa ^2 + \lpa \Hh \rpa ^2 \frac{1}{2-2H}
\le \frac{5}{2} \lpa \Hh \rpa ^2 (1-H)^{-1} .
\end{eqnarray*}
Hence taking $N=K 2^{2m}$ and $C > 1$ in (\ref{eq:basic2}), and using
(\ref{eq:logN}) too, one obtains for $m \ge 1$ that
\begin{eqnarray}
M_m &=& \max_{1 \le k \le K 2^{2m}} |A_{m,k}| \nonumber \\
&\le& \frac{5/\sqrt{2}}{\Gaa} \la \Hh \ra (1-H)^{-\half} (C \log_* K)
^\half~m^{\half} 2^{-2Hm}, \label{eq:Mmax}
\end{eqnarray}
with the exception of a set of probability at most $2 \lpa K 2^{2m} \rpa
^{1-C}$, where $K > 0$ and $C > 1$ are arbitrary.

Then substituting this and Lemma~\ref{le:main} into (\ref{eq:BHdiff}),
it follows that when $K > 0$, $C \ge 3$, and $m \ge m_4(C)$,
\begin{equation}
\max_{0 \le t \le K} |\BmpH (t) - \BmH (t)|
\le \alpha (H,K) m 2^{-\beta (H) m}
\label{eq:BHdiff2}
\end{equation}
except for an event of probability at most $8.125 (K 2^{2m})^{1-C}$
where  $\alpha (H,K)$ and $\beta (H)$ are the same as in
Lemma~\ref{le:main}. Remember that the rate of convergence in
(\ref{eq:Mmax}), just like
in parts (a) and (c) of the proof of Lemma~\ref{le:main}, is faster than
the one in parts (b) and (d) of that proof. Apart from constant
multipliers, the result of (\ref{eq:Mmax}) has the same form as the
results of (a) and (c) there. Since in the statement of this
theorem we simply replaced the faster converging factors by the slower
converging ones, the constant multipliers
of (\ref{eq:Mmax}) can be ignored if $m$ is large enough. This is why
the $\alpha(H,K)$ defined by Lemma~\ref{le:main} is suitable here too.

In the second part of the proof we compare $\BmH (t)$ to
$B_{m+j}^{(H)}(t)$, where $j \ge 1$ is an arbitrary integer. If $K > 0$,
$C \ge 3$, and $m \ge m_4(C)$, then (\ref{eq:BHdiff2}) implies that
\begin{eqnarray*}
\lefteqn{\max_{0 \le t \le K} |B_{m+j}^{(H)}(t) - \BmH (t)|
\le \sum_{k=m}^{m+j-1} \max_{0 \le t \le K}
|B_{k+1}^{(H)}(t)-B_k^{(H)}(t)|} \\
&\le& \sum_{k=m}^{\infty} \alpha (H,K) k
2^{-\beta (H) k}
\le \frac{\alpha (H,K)}{\lpa 1-2^{-\beta (H)} \rpa^2} m
2^{-\beta (H) m}.
\end{eqnarray*}
Hence one can get that
\begin{multline}
\Pb \sup_{j \ge 1} \max_{0 \le t \le K} |B_{m+j}^{(H)}(t) - \BmH
(t)| \ge \frac{\alpha (H,K)}{\lpa 1-2^{-\beta (H)} \rpa ^2} m 2^{-\beta
(H) m} \rb \\
\le \sum_{k=m}^{\infty} 8.125 \lpa K 2^{2k} \rpa ^{1-C} \le 9
\lpa K 2^{2m} \rpa ^{1-C}. \label{eq:WHdiff}
\end{multline}

By the Borel--Cantelli lemma this implies that with probability 1,
the sample paths of $\BmH (t)$ converge uniformly to a process $\WH
(t)$ on any compact interval $[0, K]$. Then $\WH (t)$ has continuous
sample paths, and inherits the properties of $\BmH (t)$ described in
Section 3: it is a centered, self-similar process with stationary
increments.  As Lemma~\ref{le:Gauss} below implies, the process $\lpa
\WH (t): t \ge 0 \rpa$ so defined is Gaussian. Therefore $\WH (t)$
is a fBM and by (\ref{eq:WHdiff}) the convergence rate of the
approximation is the one stated in the theorem. \qed
\end{proof}

The aim of the next lemma to show that integration by parts is
essentially valid for (\ref{eq:fbrow}) representing $\WH (t)$,
resulting a formula similar to (\ref{eq:fbrow2}). Then
it follows that $\lpa \WH (t): t \ge 0 \rpa$ can be
stochastically arbitrarily well approximated by a linear
transform of the Gaussian process $\lpa W(t): t \ge 0 \rpa$, so
it is also Gaussian.

\begin{lem} \label{le:Gauss}
Let $\WH (t)$ be the process whose existence is proved in Theorem~
\ref{th:main} above for $H \in (\quart, 1)$, or, by a modified
construction, in Theorem~\ref{th:mainKMT} below for any $H \in (0,
1)$. Then for any $t > 0$ and $\epsilon > 0$ there exists a
$\delta_0 > 0$ such that for any $0 < \delta < \delta_0$ we have
\begin{equation}
\Pb \la \WH (t) - \WHd (t) - \frac{\delta^{\Hh} W(t-\delta)}{\Gaa}
\ra  \ge \epsilon \rb \le \epsilon , \label{eq:WHd}
\end{equation}
where
\begin{equation}
\WHd (t) := - \int_{[-\frac{1}{\delta} , -\delta] \cup [0,
t-\delta]} h^{\prime}_{s} (s,t) \, W(s) \di s, \label{eq:WHd1}
\end{equation}
and $h(s,t)$ is defined by (\ref{eq:hst}). ($\WHd (t)$ is almost
surely well-defined pathwise as an integral of a continuous
function.)

The lemma shows that as $\delta \to 0+$, $\WHd (t)$ stochastically
converges to $\WH (t)$ when $H > \half$, while $\WHd (t)$ has a
singularity given by the extra term in (\ref{eq:WHd}) when $H <
\half$. (If $H = \half$ then $\WHd (t) = 0$ and the lemma becomes
trivial.)
\end{lem}

\begin{proof}
Fix $t > 0$ and $\epsilon > 0$ and take any $\delta$, $ 0 < \delta
\le t$. Let us introduce the notation, cf. (\ref{eq:afbrow2}):
\begin{equation}
\BmHd (\tm) = \sum_{t_r \in I_{m,\delta}} \frac{h(t_r - \Delta t, \tm)
- h(t_r, \tm)}{\Delta t} \, B_m(t_r) \, \Delta t , \label{eq:BHd}
\end{equation}
where $I_{m,\delta} = \lpa \dmi <t_r\le -\dm \rsb \cup (0 < t_r \le
\tm-\dm ]$ and the abbreviation $s_{(m)} = \lfloor s 2^{2m} \rfloor
2^{-2m}$ is used for $s=t, \delta ,$ and $-1/\delta$ (an empty sum
being zero by convention. Then we get the inequality
\begin{eqnarray}
\lefteqn{ \la \WH (t) - \WHd (t) - \frac{\delta^{\Hh}
W(t-\delta)}{\Gaa} \ra \le |\WH (t) - \BmH (\tm)|} \nonumber \\
&&+ \la \BmH (\tm) - \BmHd (\tm)
- \frac{\dm^{\Hh} B_m(\tm-\dm+\Delta t)}{\Gaa} \ra \nonumber \\
&&+ |\BmHd (\tm) - \WHd (t)| \nonumber \\
&&+ \la \frac{\dm^{\Hh} B_m(\tm-\dm+\Delta t)}{\Gaa}
- \frac{\delta^{\Hh} W(t-\delta)}{\Gaa} \ra . \label{eq:WHd2}
\end{eqnarray}

First we have to estimate the second term on the right hand side as
$\delta \to 0+$, uniformly in $m$ (this requires the longest
computation):
\[ \BmH (\tm) - \BmHd (\tm) - \frac{\dm^{\Hh}
B_m(\tm-\dm+\Delta t)}{\Gaa} =: E_{m,\delta} + F_{m,\delta}
+ G_{m,\delta}, \]
where
\begin{eqnarray*}
E_{m,\delta} &=&  \sum_{\tm - \dm < t_r \le \tm} \frac{h(t_r
- \Delta t, \tm) - h(t_r, \tm)}{\Delta t} \, B_m(t_r) \, \Delta t \\
&& \qquad - \frac{\dm^{\Hh} B_m(\tm-\dm+\Delta t)}{\Gaa},
\end{eqnarray*}
\[ F_{m,\delta} = \sum_{-\dm < t_r \le 0} \frac{h(t_r
- \Delta t, \tm) - h(t_r, \tm)}{\Delta t} \, B_m(t_r) \, \Delta t
\]
and
\[ G_{m,\delta} = \sum_{-\infty < t_r \le \dmi}
\frac{h(t_r - \Delta t, \tm) - h(t_r, \tm)}{\Delta t} \,
B_m(t_r) \, \Delta t . \]

Then ``summation by parts'' shows that
\[ E_{m,\delta} = \sum_{\tm - \dm < t_r < \tm} h(t_r, \tm )
[B_m(t_{r+1}) - B_m(t_r)] . \]
(This is the point where the extra term in the definition
of $E_{m,\delta}$ is needed.)
Thus
\begin{eqnarray*}
\lefteqn{\Var \lpa \Gaa E_{m,\delta} \rpa = \sum_{\tm - \dm < t_r < \tm}
(\tm - t_r)^{2H-1} 2^{-2m}} \\
&=& \tm^{2H-1} \sum_{\tm - \dm < t_r < \tm}
\lpa 1 - \frac{t_r}{\tm}\rpa^{2H-1} \Delta t \\
&\le& \tm^{2H-1} \int_{\tm - \dm}^{\tm} \lpa 1 - \frac{u}{\tm} \rpa \di u \\
&=& \frac{\dm^{2H}}{2H} \le \frac{\delta^{2H}}{2H},
\end{eqnarray*}
for any $m \ge 0$.
Then by the large deviation inequality (\ref{eq:ldp}), for any $m \ge 0$
and for any $C>0$,
\begin{equation}
\Pb |E_{m,\delta}| \ge \lpa 2C \log_* \lpa \frac{1}{\delta} \rpa
\rpa^{\half} \frac{\delta^H}{\Gaa (2H)^{\half}} \rb \le 2 \delta^C .
\label{eq:Emd}
\end{equation}

Similarly as above, the definition of $F_{m,\delta}$ can be rewritten
using ``summation by parts'' that gives
\[ F_{m,\delta} = \sum_{-\dm < t_r < 0} h(t_r, \tm )
[B_m(t_{r+1}) - B_m(t_r)] + h(-\dm , \tm ) B_m(-\dm + \Delta t). \]
The definition of $F_{m,\delta}$ shows that it is equal to zero whenever
$\delta < \Delta t$, therefore when giving an upper bound for its
variance it can be assumed that $\delta \ge \Delta t$. Thus
\begin{eqnarray*}
\lefteqn{\Var \lpa \Gaa F_{m,\delta} \rpa} \\
&=& \sum_{0<t_v<\dm} \lsb (\tm+t_v)^{\Hh} - t_v^{\Hh} \rsb^2 \Delta t
 + \lsb (\tm+\dm)^{\Hh} - \dm^{\Hh} \rsb^2 \\
 &&\times(\dm - \Delta t) \\
&\le& \tm^{2H-1} \sum_{0<t_v<\dm} \lsb \lpa 1+\frac{t_v}{\tm}\rpa^{2H-1}
+ \lpa \frac{t_v}{\tm} \rpa^{2H-1} \rsb \Delta t \\
&&+ \lsb (\tm+\dm)^{2H-1} + \dm^{2H-1} \rsb \dm \\
&\le& \tm^{2H-1} \int_0^{\dm} \lsb \lpa 1+\frac{u}{\tm}\rpa^{2H-1}
+ \lpa \frac{u}{\tm} \rpa^{2H-1} \rsb \, dt + 2 \tm^{2H-1} \dm +
\dm^{2H} \\
&=& \frac{\tm^{2H}}{2H} \lsb \lpa 1+\frac{\dm}{\tm}\rpa^{2H}
+ \lpa \frac{\dm}{\tm} \rpa^{2H} - 1 \rsb + 2 \tm^{2H-1} \dm + \dm^{2H}
\\
&\le& \threehalf \tm^{2H-1} \dm + \frac{\dm^{2H}}{2H} + 2 \tm^{2H-1} \dm
+ \dm^{2H} \le \frac72 t^{2H-1} \delta +\frac{3}{2H} \delta^{2H} .
\end{eqnarray*}
So by the large deviation inequality (\ref{eq:ldp}), for any $m \ge 0$
and for any $C>0$,
\begin{equation}
\Pb |F_{m,\delta}| \ge \lpa 2C \log_* \lpa \frac{1}{\delta} \rpa
\rpa^{\half} \frac{\lpa \frac72 t^{2H-1} \delta +\frac{3}{2H} \delta^{2H}
\rpa^{\half}}{\Gaa} \rb \le 2 \delta^C . \label{eq:Fmd}
\end{equation}

Proceeding in a similar way with $G_{m,\delta}$, one obtains that
\begin{eqnarray*}
G_{m,\delta} &=& \sum_{-\infty < t_r < \dmi} h(t_r, \tm )
[B_m(t_{r+1}) - B_m(t_r)] \\
&&- h\lpa \dmit,\tm \rpa B_m\lpa \dmit \rpa .
\end{eqnarray*}
Hence
\begin{eqnarray*}
\lefteqn{\Var (\Gaa G_{m,\delta}) = \sum_{-\dmi<t_v<\infty} \lsb
(\tm+t_v)^{\Hh} - t_v^{\Hh} \rsb^2 \Delta t} \\
&& + \lsb \lpa \tm-\dmit \rpa^{\Hh}
- \lpa -\dmit \rpa^{\Hh} \rsb^2 \lpa -\dmit \rpa \\
&\le& \sum_{-\dmi<t_v<\infty} t_v^{2H-1} \lpa H-\half \rpa^2 \lpa
\frac{\tm}{t_v} \rpa^2 \Delta t \\
&& + \lpa -\dmit \rpa^{2H-1} \lpa \Hh
\rpa^2 \lpa \frac{\tm}{-\dmit} \rpa^2 \lpa -\dmit \rpa \\
&=& \lpa \Hh \rpa^2 \tm^2 \lb \sum_{-\dmi<t_v<\infty} t_v^{2H-3} \Delta t
+ \lpa -\dmit \rpa^{2H-2} \rb \\
&\le& \lpa \Hh \rpa^2 \tm^2 \lb \int_{-\dmi}^{\infty} u^{2H-3} \di u
+ \lpa -\dmit \rpa^{2H-2} \rb \\
&\le& \frac{3 (\Hh)^2}{2(1-H)} t^2 \delta^{2-2H}.
\end{eqnarray*}
So again by the large deviation inequality (\ref{eq:ldp}),
for any $m \ge 0$ and for any $C>0$,
\begin{equation}
\Pb |G_{m,\delta}| \ge \lpa 2C \log_* \lpa \frac{1}{\delta} \rpa
\rpa^{\half} \frac{\la \Hh \ra}{\Gaa} \lpa \frac{\threehalf}{1-H} \rpa^{\half}
t \,\delta^{1-H} \rb \le 2 \delta^C . \label{eq:Gmd}
\end{equation}

Combining (\ref{eq:Emd}), (\ref{eq:Fmd}) and (\ref{eq:Gmd}), it follows
that there exists a $\delta_0 > 0$ such that for any $0<\delta <\delta_0$
and for any $m \ge 0$,
\[ \Pb \la \BmH (\tm) - \BmHd (\tm) - \frac{\dm^{\Hh} B_m(\tm-\dm
+\Delta t)}{\Gaa} \ra \ge \frac{\epsilon}{4} \rb \le \frac{\epsilon}{4} .
\]

After the second term on the right hand side of (\ref{eq:WHd2}) we
turn to the third term. Take now any $\delta \in (0, \delta_0)$
Since $h(s,t)$ has continuous partial derivative w.r.t. $s$ on
the intervals $[-1/\delta, -\delta]$ and $[\delta, t-\delta]$ and
by Theorem~\ref{th:wiener}, $B_m$ a.s. uniformly converges to the
Wiener process $W$ on these intervals, comparing (\ref{eq:WHd1}) and
(\ref{eq:BHd}) shows that with this $\delta$ there exists an $m$ such
that
\[ \Pb |\BmHd (\tm) - \WHd (t) | \ge \frac{\epsilon}{4} \rb \le
\frac{\epsilon}{4} . \]
Theorem~\ref{th:wiener} also implies that $m$ can be chosen so that
for the fourth term in (\ref{eq:WHd2}) one similarly has
\[ \Pb \la \frac{\dm^{\Hh} B_m(\tm-\dm+\Delta t)}{\Gaa}
- \frac{\delta^{\Hh} W(t-\delta)}{\Gaa} \ra \ge \frac{\epsilon}{4}
\rb \le \frac{\epsilon}{4} . \]
Finally, Theorem~\ref{th:main} (or, with a modified construction,
Theorem~\ref{th:mainKMT} below) guarantees that $m$ can be chosen so
that the first term in (\ref{eq:WHd2}) satisfies the same inequality:
\[ \Pb |\WH (t) - \BmH (t)| \ge \frac{\epsilon}{4} \rb \le
\frac{\epsilon}{4} . \]
The last four formulae together prove the lemma. \qed
\end{proof}


\section{Improved construction using the KMT approximation}

Parts (b) and (d) of the proof of Lemma~\ref{le:main} gave worse rate of
convergence than parts (a) and (c), in which the rates can be conjectured
to be best possible. The reason for this is clearly the relatively weaker
convergence rate of the RW approximation of ordinary BM, that was used in
parts (b) and (d), but not in parts (a) and (c). It is also clear from
there that using the best possible KMT approximation instead would
eliminate this weakness and would give hopefully the best possible rate
here too. The price one has to pay for this is the intricate
and ``future-dependent'' procedure by which the KMT method
constructs suitable approximating RWs from BM.

The result we need from Koml\'os, Major, and Tusn\'ady (1975, 1976) is
as follows. Suppose that one wants to define an i.i.d. sequence $X_1,
X_2, \dots$ of random variables with a given distribution so that the
partial sums are as close to BM as possible. Assume that
$\E (X_k)=0$, $\Var (X_k)=1$ and the moment generating function $\E \lpa
e^{u X_k} \rpa < \infty$ for $|u| \le u_0, \ u_0 > 0$.  Let $S(k) = X_1+
\cdots +X_k$, $k \ge 1$ be the partial sums. If BM $W(t)$ $(t \ge 0)$ is
given, then for any $n \ge 1$ there exists a sequence of conditional
quantile transformations applied to $W(1), W(2), \dots , W(n)$ so that
one obtains the desired partial sums $S(1), S(2), \dots , S(n)$ and the
difference between the two sequences is the smallest possible:
\begin{equation}
\Pb \max_{0 \le k \le n} |S(k) - W(k)| > C_0 \log n + x \rb <
K_0 e^{-\lambda x}, \label{eq:KMT0}
\end{equation}
for any $x>0$, where $C_0, K_0, \lambda$ are positive
constants that may depend on the distribution of $X_k$, but not on
$n$ or $x$. Moreover, $\lambda$ can be made arbitrarily large by
choosing a large enough $C_0$. Taking $x=C_0 \log n$ here one obtains
\begin{equation}
\Pb \max_{0 \le k \le n} |S(k) - W(k)| > 2 C_0 \log n \rb <
K_0 n^{-\lambda C_0}, \label{eq:KMT}
\end{equation}
where $n \ge 1$ is arbitrary.

Fix an integer $m \ge 0$, and introduce the same notations
as in previous sections: $\Delta t = 2^{-2m}$, $t_x=x \Delta t$. Then
multiply the inner inequality in (\ref{eq:KMT}) by $2^{-m}$ and use
self-similarity (\ref{eq:ss}) of BM (with $H=\half$) to obtain a shrunken
RW $B_m^*(t_k)=2^{-m}S_m(k)$ $(0 \le k \le K2^{2m})$ from the
corresponding dyadic values $W(t_k)$ $(0 \le k \le K2^{2m})$ of BM by a
sequence of conditional quantile transformations so that
\begin{equation}
\max_{0 \le t_k \le K} |B_m^*(t_k) - W(t_k)| \le 2C_0 2^{-m}
\log_* (K 2^{2m}) \le 5C_0 \log_* K \: m 2^{-m} , \label{eq:KMT1}
\end{equation}
with the exception of a set of probability smaller than $K_0 (K 2^{2m})
^{-\lambda C_0}$, for any $m \ge 1$ and $K>0$. [Here (\ref{eq:logN})
was used too.] Then (\ref{eq:KMT1}) implies for the difference of two
consecutive approximations that
\begin{multline}
\Pb \max_{0 \le t_k \le K} |B_{m+1}^*(t_k) - B_m^*(t_k)|
> 10C_0 \log_* K \: m 2^{-m} \rb \\
< 2K_0 \lpa K 2^{2m} \rpa ^{-\lambda C_0}  \label{eq:KMT2}
\end{multline}
for any $m \ge 1$ and $K>0$. This is exactly that we need to improve
the rates of convergence in parts (b) and (d) of Lemma~\ref{le:main}.

Substitute these KMT approximations $B_m^*(t_r)$ into definition
(\ref{eq:afbrow1}) or (\ref{eq:afbrow2}) of $\BmH (t_k)$. This way one
can obtain faster converging approximations $B_m^{*(H)}$ of fBM. Then
everything above in Sections 3 and 4 are still valid, except that one
can use the improved formula (\ref{eq:KMT2}) instead of
Lemma~\ref{le:apprb} at parts (b) and (d) in the proof of
Lemma~\ref{le:main}. This way, instead of (\ref{eq:Ymax}) one gets
\begin{equation}
\max_{1 \le k \le K 2^{2m}} |Y_{m,k}| \le \left\{
\begin{array}{ll}
23 C_0 \log_* K \: m 2^{-2Hm}& \mbox{ if } 0 < H < \half , \\
15 C_0 \log_* K \: K^{\Hh} m 2^{-m}&
\mbox{ if } \half < H < 1,
\end{array} \right. \label{eq:YmKMT}
\end{equation}
for any $m \ge 1$, except for a set of probability smaller than
$2K_0 (K 2^{2m})^{-\lambda C_0}$.

Also by (\ref{eq:KMT2}), instead of (\ref{eq:Ufact1}) and (\ref{eq:Uprob})
one has the improved inequalities:
\begin{eqnarray}
\lefteqn{\max_{(j-1)L<t_v\le jL} \la B_{m+1}^*(-t_v) - B_m^*(-t_v) \ra}
\nonumber \\
&\le& \left\{
\begin{array}{ll}
10 C_0 \log_*L\: m 2^{-m} & \mbox{ if } j=1,  \\
14 C_0 \log_*j\: \log_*L\: m 2^{-m} & \mbox{ if } j \ge 2,
\end{array} \right. \label{eq:Uf1KMT}
\end{eqnarray}
with the exception of a set of probability smaller than $2K_0 (jL 2^{2m})
^{-\lambda C_0}$, where $m \ge 1$. If $C_0$ is chosen large enough
so that $\lambda C_0 \ge 2$, then (\ref{eq:Uf1KMT}) holds simultaneously
for all $j=1,2,3, \dots$  except for a set of probability smaller than
\begin{equation}
2K_0 (L 2^{2m})^{-\lambda C_0} \sum_{j=1}^{\infty} j^{-\lambda C_0}
< \quart K_0 (K 2^{2m})^{-\lambda C_0}. \label{eq:UpKMT}
\end{equation}
(Remember that we chose $L=4K$ in part (d) of the proof of
Lemma~\ref{le:main}.) Then using this in part (d) of Lemma~\ref{le:main},
instead of (\ref{eq:Uj}) one needs the estimate
\[ \sum_{j=2}^{\infty} j^{H-\frac52} \log_* j < \int_1^{\infty}
x^{H-\frac52} \log_* x \di x < 1 \qquad (0<H<1). \]
Then instead of (\ref{eq:Umax1}) and (\ref{eq:Umax2}), the improved
results are as follows. First, in the case $\half <H<1$ one has
\begin{eqnarray}
\lefteqn{\max_{1 \le k \le K 2^{2m}} |U_{m,k}|} \nonumber \\
&\le& 10 C_0 \log_*L\: m 2^{-m} \frac{3}{8} L^{\Hh} + 14 C_0
\log_*L\: m 2^{-m} \frac{5}{3}~ \la \Hh \ra L^{\Hh} \nonumber \\
&\le& \lpa 18 + 112 ~ \la \Hh \ra \rpa C_0 \log_*K\: K^{H-\half}
m 2^{-m} \label{eq:Um1KMT}
\end{eqnarray}
for any $m \ge 1$ and $C_0$ large enough so that $\lambda C_0 \ge 2$,
except for a set of probability smaller than given by (\ref{eq:UpKMT}).

Now in the case $0<H<\half$ it follows that
\begin{eqnarray}
\lefteqn{\max_{1 \le k \le K 2^{2m}} |U_{m,k}|} \nonumber \\
&\le& 10 C_0 \log_*L\: m 2^{-m} \threehalf(\Delta t)^{\Hh}
+ 14 C_0 \log_*L\: m 2^{-m} \, \frac{5}{2}~\la \Hh \ra L^{\Hh}
\nonumber \\
&\le& C_0 \log_*K\: m \lpa 36 \cdot 2^{-2Hm} + 84 \, \la \Hh \ra
K^{H-\half} 2^{-m} \rpa \label{eq:Um2KMT}
\end{eqnarray}
for any $m \ge 1$ and $C_0$ large enough so that $\lambda C_0 \ge 2$,
except for a set of probability smaller than given by (\ref{eq:UpKMT}).

As a result, there is convergence for any $H \in (0, 1)$. Since the
KMT approximation itself has best possible rate for approximating
ordinary BM by RW, it can be conjectured that the resulting convergence
rates in the next lemma and theorem are also best possible (apart from
constant multipliers) for approximating fBM by moving averages of a RW.

\begin{lem} \label{le:mainKMT}
For any $H \in (0, 1)$, $m \ge 1$, $K > 0$, $C > 1$, and $C_0$ large
enough, we have
\begin{eqnarray*}
\lefteqn{\Pb \max_{0 \le t_k \le K} |B_{m+1}^{*(H)} (t_k) -
B_m^{*(H)} (t_k)| \ge \alpha^* m 2^{-\beta^*(H) m} \rb} \\
&\le& 4 (K 2^{2m})^{1-C} + 3 K_0 (K 2^{2m})^{-\lambda C_0},
\end{eqnarray*}
where $t_k = k 2^{-2m}$, $\beta^*(H)
= \min (2H, 1)$, $\alpha^* = \alpha^*(H,K,C,C_0)$,
\[
\alpha^*
= \frac{(\log_* K)^{\half}}{\Gaa} \lsb 10 C^{\half} \frac{\la \Hh \ra}
{(1-H)^{\half}} + C_0(\log_* K)^{\half} \lpa 59 + 84 \la \Hh \ra
K^{H-\half} \rpa \rsb
\]
if $H \in \lpa 0 , \half \rpa$ ,
\[
\alpha^*
= \frac{(\log_* K)^{\half}}{\Gaa} \lsb 10 C^{\half} \frac{ \la \Hh \ra}
{(1-H)^{\half}} + C_0 (\log_* K)^{\half} \lpa 33 + 112 \la \Hh \ra \rpa
K^{H-\half} \rsb
\]
if $H \in \lpa \half , 1 \rpa$, and the constants $\lambda$,
$C_0$ and $K_0$ are defined by the KMT approximation (\ref{eq:KMT0}) with
$C_0$ is chosen so large that $\lambda C_0 \ge 2$.
[The case $H=\half$ is described by (\ref{eq:KMT2}).]
\end{lem}

\begin{proof}
Combine the results of parts (a) and (c) in the proof of Lemma~
\ref{le:main} and the improved inequalities above, that is, apply
(\ref{eq:ZYVU}), (\ref{eq:Zmax}), (\ref{eq:YmKMT}),
(\ref{eq:Vmax}), (\ref{eq:Um1KMT}), and (\ref{eq:Um2KMT}). Here too,
we simply replace the faster converging factors by the slower
converging ones, but the constant multipliers of faster converging terms
cannot be ignored, since the lemma is stated for any $m \ge 1$. \qed
\end{proof}

Now we can extend the improved approximations of fBM to real arguments
by linear interpolation, in the same way as we did with the original
approximations, see (\ref{eq:intpol}). This way we
get continuous parameter approximations $B_m^{*(H)}(t)$, $(t \ge 0)$ for
$m=0,1,2, \dots$, with continuous, piecewise linear sample paths.
Now we can state the second main result of this paper.

\begin{thmn} \label{th:mainKMT}
For any $H \in (0, 1)$, the sequence $B_m^{*(H)}(t)$  $(t \ge 0, m=0,1,2,
\dots )$ a.s. uniformly converges to a fBM $\WH (t)$  $(t \ge 0)$ on any
compact interval $[0, K]$, $K > 0$. If $m \ge 1$, $K > 0$, $C \ge 2$,
and $C_0$ is large enough, it follows that
\begin{eqnarray*}
\lefteqn{\Pb \max_{0 \le t \le K} |\WH (t) - B_m^{*(H)}(t)|
\ge \frac{\bar{\alpha}^*}{\lpa 1-2^{-\beta^*(H)} \rpa ^2} m
2^{-\beta^*(H) m} \rb} \\
&& \qquad \le 6 (K 2^{2m})^{1-C} + 4 K_0 (K 2^{2m})^{-\lambda C_0}
\end{eqnarray*}
with
\[ \bar{\alpha}^* = \alpha^* + \frac{10}{\Gaa} C^{\half} (\log_*K)
^{\half} \la \Hh \ra (1-H)^{-\half}, \]
where $\alpha^*$ and $\beta^*(H)$ are the same as in Lemma~\ref
{le:mainKMT}. (In other words, in the definition of $\alpha^*$ in
Lemma~\ref{le:mainKMT} the constant multiplier 10 has to be changed
to 20 here.) The constants $\lambda, C_0, K_0$ are defined by the KMT
approximation (\ref{eq:KMT0}) with $C_0$ is chosen so large that
$\lambda C_0 \ge 2$. [The case $H=\half$ is described by
(\ref{eq:KMT1}).]
\end{thmn}

\begin{proof}
The proof can follow the line of the proof of Theorem~\ref{th:main}
with one exception: the constant multipliers in (\ref{eq:Mmax}) and
consequently in (\ref{eq:BHdiff}) cannot be ignored here.
This is why the multiplier $\alpha^*$ of Lemma~\ref{le:mainKMT}
had to be modified in the statement of the theorem. \qed
\end{proof}

It can be conjectured that the best rate of approximation of fBM by
moving averages of simple RW's is $O(N^{-H} \log N)$, where $N$ is the
number of points considered. Though it seems quite possible that
definition of $B_m^{*(H)}(t)$ above, see (\ref{eq:afbrow1}) with the KMT
approximations $B_m^*(t_r)$, supplies this rate of convergence for any
$H \in (0, 1)$, but in Theorem~\ref{th:mainKMT} we were able to prove
this rate only when $H \in (0, \half)$. A possible explanation could be
that in parts (b) and (d) of Lemma~\ref{le:main} we separated the maxima
of the kernel and the ``integrator'' parts.

As a result, the convergence rate we were able to prove when $\half <
H < 1$ is the same $O(N^{-\half} \log N)$ that the original KMT
approximation (\ref{eq:KMT1}) gives for ordinary BM, where $N=K
2^{2m}$, though in this case the sample paths of fBM are smoother
than that of BM. (See e.g. \cite{Dec98}.) On the other hand, the
obtained convergence rate is worse than this, but still thought to be
the best possible, $O(N^{-H} \log N)$, when $0 < H < \half$, which
heuristically can be explained by the more zigzagged sample paths of
fBM in this case.



\end{document}